# ON THE ROLE OF QUADRATIC OSCILLATIONS IN NONLINEAR SCHRÖDINGER EQUATIONS

RÉMI CARLES, CLOTILDE FERMANIAN KAMMERER, AND ISABELLE GALLAGHER

ABSTRACT. We consider a nonlinear semi–classical Schrödinger equation for which it is known that quadratic oscillations lead to focusing at one point, described by a nonlinear scattering operator. If the initial data is an energy bounded sequence, we prove that the nonlinear term has an effect at leading order *only if* the initial data have quadratic oscillations; the proof relies on a linearizability condition (which can be expressed in terms of Wigner measures). When the initial data is a sum of such quadratic oscillations, we prove that the associate solution is the superposition of the nonlinear evolution of each of them, up to a small remainder term. In an appendix, we transpose those results to the case of the nonlinear Schrödinger equation with harmonic potential.

## 1. INTRODUCTION

Consider the initial value problem

(1.1) $$\begin{cases} i\partial_t u + \frac{1}{2}\Delta u = 0, & (t,x) \in \mathbb{R}_+ \times \mathbb{R}^n, \\ u_{|t=0} = e^{-i\frac{|x|^2}{2}}. \end{cases}$$

It is easy to see that at time $t=1$, the solution $u$ is the Dirac mass at the origin. Robbiano and Zuily ([RZ00]) proved that the presence of quadratic oscillations is essentially the only cause in the formation of singularities in linear Schrödinger equations. They considered the equation

$$\begin{cases} i\partial_t u + \Delta_g u = 0, \\ u_{|t=0} = f(x)e^{i\varphi(x)}, \end{cases}$$

where $\Delta_g$ is the Laplacian on an asymptotically flat analytic metric. The amplitude $f$ is in $L^2(\mathbb{R}^n)$, and the phase $\varphi$ is real analytic. If $\varphi = P_m + R_m$, where $P_m$ is a homogeneous polynomial of degree $m$, and $R_m(z) = o(|z|^m)$ as $|z| \to \infty$, then analytic singularities can appear for positive times only in the following two cases: $m=2$, $P_2(y) = \langle Ay, y\rangle$ and $A$ has a real negative eigenvalue, or $m \geq 3$ and $\partial_y P_m$ has a real zero. The case $m=2$ shows that the example (1.1) is in some sense generic.

For the nonlinear Schrödinger equation with attractive nonlinearity,

$$\begin{cases} i\partial_t u + \frac{1}{2}\Delta u = -|u|^{2\sigma}u, & (t,x) \in \mathbb{R}_+ \times \mathbb{R}^n, \\ u_{|t=0} = u_0(x), \end{cases}$$

where $\sigma \in [2/n, 2/(n-2)[$, and $u_0 \in H^1(\mathbb{R}^n)$ with $|x|u_0 \in L^2(\mathbb{R}^n)$, it is well known (cite e.g. [Caz93]) that blow up in finite time may occur. Cazenave and Weissler ([CW92]) proved that changing $u_0(x)$ into $u_0(x)e^{ib|x|^2}$ for $b$ sufficiently large ensures the existence of the solution $u$ globally in time. On the other hand, if $u$ blows up

---

2000 *Mathematics Subject Classification*. Primary 35Q55; Secondary 35B40, 35B05.

This work was partially supported by the ACI grant "Équation des ondes : oscillations, dispersion et contrôle".

in finite time because it has negative energy, then changing $u_0(x)$ into $u_0(x)e^{-ib|x|^2}$ for $b>0$ sufficiently large makes the blow up happen sooner.

For the Schrödinger equation with focusing nonlinearity and critical power,

$$(1.2) \qquad i\partial_t u + \frac{1}{2}\Delta u = -|u|^{4/n}u,$$

Merle ([Mer93]) proved that if the initial data $u_0$ belongs to $H^1(\mathbb{R}^n)$ and its $L^2$-norm is the same as that of the solitary wave, then blow up in finite time $T>0$ can occur only under a very rigid assumption on $u_0$, for it must be of the form

$$u_0(x) = \left(\frac{\delta}{T}\right)^{n/2} e^{i\theta - i|x-x_1|^2/2T + i\delta^2/T} Q\left(\delta\left(\frac{x-x_1}{T} - x_0\right)\right),$$

for some $\theta \in \mathbb{R}$, $\delta > 0$, $x_0, x_1 \in \mathbb{R}^n$, where $Q$ is a solution of the stationary problem. In particular, blow up in finite time necessarily involves quadratic oscillations in the initial data (but this condition is not sufficient).

In [MV98], the authors notice that the defect of compactness in the two-dimensional cubic nonlinear Schrödinger equation is due to quadratic oscillations, $e^{-i\lambda x^2}$, with $\lambda$ large. These oscillations are related to the Galilean invariance, as in [Mer93].

In [Car00], nonlinear Schrödinger equations are considered in the semiclassical limit, in particular the initial value problem,

$$(1.3) \qquad \begin{cases} i\varepsilon \partial_t u^\varepsilon + \frac{1}{2}\varepsilon^2 \Delta u^\varepsilon = \varepsilon^{n\sigma}|u^\varepsilon|^{2\sigma} u^\varepsilon, & (t,x) \in \mathbb{R}_+ \times \mathbb{R}^n, \\ u^\varepsilon_{|t=0} = e^{-i\frac{|x|^2}{2\varepsilon}} f(x), \end{cases}$$

where $\varepsilon \in ]0,1]$, $\sigma > 2/(n+2)$, and $\sigma < 2/(n-2)$ if $n \geq 3$. The asymptotic behavior of the solution $u^\varepsilon$ is studied, as $\varepsilon$ goes to zero: quadratic oscillations cause focusing at the origin at time $t=1$ (compare with (1.1)), and the scaling of the nonlinearity (in particular, the presence of the factor $\varepsilon^{n\sigma}$) makes the influence of the right-hand side of (1.3) negligible away from the focal point. On the other hand, the caustic crossing takes some nonlinear effects into account, and is described at leading order by the (nonlinear) scattering operator associated with

$$(1.4) \qquad i\partial_t \psi + \frac{1}{2}\Delta \psi = |\psi|^{2\sigma}\psi.$$

One may argue that this case is very particular, inasmuch as the initial oscillations are associated with a specific geometry (the rays of geometric optics meet at the origin at time 1), and ask: What if the initial data are more general? For instance, what happens when more general oscillating initial data are considered,

$$(1.5) \qquad u^\varepsilon_{|t=0} = \sum_{j=1}^J f_j(x) e^{i\frac{\varphi_j(x)}{\varepsilon}}?$$

The aim of this paper is to study (a generalization of) these cases. We prove in particular that the nonlinear term $\varepsilon^{n\sigma}|u^\varepsilon|^{2\sigma}u^\varepsilon$ can always be neglected if none of the phases $\varphi_j$ is quadratic (Theorem 1.2): this means that the framework of [Car00] concerned a critical case as far as geometric optics is concerned. We also study the case where all the phases $\varphi_j$ are quadratic (Theorem 1.4), and prove that the solution $u^\varepsilon$ can be described as the superposition of the solutions $v^\varepsilon_j$ of problems (1.3), with

$$v^\varepsilon_{j|t=0} = f_j(x) e^{i\frac{\varphi_j(x)}{\varepsilon}}.$$

We now set up the framework we will keep throughout this paper, and state precisely our results. We study the following Cauchy problem,

$$
(1.6) \quad \begin{cases} i\varepsilon \partial_t u^\varepsilon + \frac{1}{2}\varepsilon^2 \Delta u^\varepsilon = \varepsilon^{n\sigma}|u^\varepsilon|^{2\sigma} u^\varepsilon, & (t,x) \in \mathbb{R}_+ \times \mathbb{R}^n, \\ u^\varepsilon_{|t=0} = u^\varepsilon_0, \end{cases}
$$

with the following assumptions.

**Assumptions.** *We suppose that*

- $(H1)$   $\varepsilon \in ]0,1]$.
- $(H2)$   $\sigma > 2/n$, and $\sigma < 2/(n-2)$ if $n \geq 3$.
- $(H3)$   *The initial data $u^\varepsilon_0$ belong to $H^1(\mathbb{R}^n)$, uniformly in the following sense,*

$$\sup_{0<\varepsilon\leq 1}\left(\|u^\varepsilon_0\|_{L^2} + \|\varepsilon \nabla u^\varepsilon_0\|_{L^2}\right) < \infty.$$

  *Notice that this case includes the case of WKB data (1.5).*
- $(H4)$   *There is no focusing at time 0. As we shall see later, that means that we suppose that*

$$\limsup_{\varepsilon \to 0} \varepsilon^{n\sigma} \|u^\varepsilon_0\|^{2\sigma+2}_{L^{2\sigma+2}} = 0.$$

We define $v^\varepsilon$ as the free evolution of $u^\varepsilon_0$,

$$
(1.7) \quad \begin{cases} i\varepsilon \partial_t v^\varepsilon + \frac{1}{2}\varepsilon^2 \Delta v^\varepsilon = 0, & (t,x) \in \mathbb{R}_+ \times \mathbb{R}^n, \\ v^\varepsilon_{|t=0} = u^\varepsilon_0. \end{cases}
$$

We will also need the free evolution "without $\varepsilon$",

$$
(1.8) \quad \begin{cases} i\partial_t V + \frac{1}{2}\Delta V = 0, & (t,x) \in \mathbb{R}_+ \times \mathbb{R}^n, \\ V_{|t=0} = V_0. \end{cases}
$$

We will define the associate linear operators $U^\varepsilon_0(t) := e^{i\varepsilon\frac{t}{2}\Delta}$ and $U_0(t) := e^{i\frac{t}{2}\Delta}$.

Following some ideas introduced in [Gér96], our first result is the following.

**Theorem 1.1.** *Let $T > 0$. The following properties are equivalent,*
*(1) The function $v^\varepsilon$ is an approximation of $u^\varepsilon$ on the time interval $[0,T]$,*

$$\sup_{0\leq t\leq T}\left(\|u^\varepsilon(t) - v^\varepsilon(t)\|_{L^2} + \|\varepsilon \nabla u^\varepsilon(t) - \varepsilon \nabla v^\varepsilon(t)\|_{L^2}\right) \xrightarrow[\varepsilon \to 0]{} 0.$$

*(2) The function $v^\varepsilon$ satisfies*

$$(1.9) \quad \limsup_{\varepsilon \to 0} \sup_{0\leq t\leq T} \varepsilon^{n\sigma} \|v^\varepsilon(t)\|^{2\sigma+2}_{L^{2\sigma+2}} = 0.$$

To check whether the condition (1.9) is satisfied or not, one can compute the Wigner measure of the initial data $u^\varepsilon_0$. We shall state the corresponding result after having analyzed more precisely condition (1.9) in Theorem 1.2.

For the sake of readability, we introduce the following notation.

**Notation.** i) For a family $(a^\varepsilon)_{0<\varepsilon\leq 1}$ of functions in $H^1(\mathbb{R}^n)$, define

$$\|a^\varepsilon\|_{H^1_\varepsilon} := \|a^\varepsilon\|_{L^2} + \|\varepsilon \nabla a^\varepsilon\|_{L^2},$$

$$\|a^\varepsilon\|_{L^{2\sigma+2}_\varepsilon} := \varepsilon^{\frac{n\sigma}{2\sigma+2}}\|a^\varepsilon\|_{L^{2\sigma+2}}.$$

We will say that $a^\varepsilon$ is bounded (resp. goes to zero) in $H^1_\varepsilon$ if

$$\limsup_{\varepsilon \to 0}\|a^\varepsilon\|_{H^1_\varepsilon} < \infty \text{ (resp. } = 0\text{)}.$$

The same notions are obviously defined in $L_\varepsilon^{2\sigma+2}$.

ii) If $(\alpha^\varepsilon)_{0<\varepsilon\leq 1}$ and $(\beta^\varepsilon)_{0<\varepsilon\leq 1}$ are two families of positive numbers, we will write

$$\alpha^\varepsilon \lesssim \beta^\varepsilon \tag{1.10}$$

if there exists $C$ independent of $\varepsilon \in ]0,1]$ (but possibly depending on other parameters) such that for any $\varepsilon \in ]0,1]$,

$$\alpha^\varepsilon \leq C\beta^\varepsilon.$$

*Remark.* Notice that (1.9) means exactly that $v^\varepsilon$ goes to zero in $L_\varepsilon^{2\sigma+2}$, uniformly for $t \in [0,T]$.

Now we examine the case where (1.9) is not satisfied, that is when the evolution of $u^\varepsilon$ on $[0,T]$ takes some nonlinear effects into account at leading order. To state our result, we adapt some techniques developed in [Gér98], [BG99] and [Ker01] for instance. Before stating the result, let us give the following definition: if $(z_j^\varepsilon)_{j\in\mathbb{N}}$ is a family of sequences in $\mathbb{R}^d$, $d \geq 1$, then we shall say that it is an orthogonal family if

$$\forall j \neq k, \quad \limsup_{\varepsilon \to 0} \frac{|z_j^\varepsilon - z_k^\varepsilon|}{\varepsilon} = \infty. \tag{1.11}$$

**Theorem 1.2.** *Let $T > 0$ and assume that (1.9) is not satisfied. Then up to the extraction of a subsequence, there exist an orthogonal family $(t_j^\varepsilon, x_j^\varepsilon)_{j\in\mathbb{N}}$ in $\mathbb{R}_+ \times \mathbb{R}^n$, a family $(\Psi_\ell^\varepsilon)_{\ell\in\mathbb{N}}$, bounded in $H_\varepsilon^1(\mathbb{R}^n)$, and a (nonempty) family $(\varphi_j)_{j\in\mathbb{N}}$, bounded in $\{\varphi \in L^2(\mathbb{R}^n) \,;\, |x|\varphi \in L^2(\mathbb{R}^n)\}$, such that:*

$$u_0^\varepsilon(x) = \Psi_\ell^\varepsilon(x) + r_\ell^\varepsilon(x), \quad \text{with } \limsup_{\varepsilon \to 0} \|U_0^\varepsilon(t) r_\ell^\varepsilon\|_{L^\infty(\mathbb{R}_+, L_\varepsilon^{2\sigma+2})} \xrightarrow[\ell\to\infty]{} 0, \tag{1.12}$$

*and for every $\ell \in \mathbb{N}$, the following asymptotics holds in $L^2(\mathbb{R}^n)$, as $\varepsilon \to 0$,*

$$\Psi_\ell^\varepsilon(x) = \sum_{j=0}^\ell \frac{1}{(t_j^\varepsilon)^{\frac{n}{2}}} \varphi_j\left(\frac{x - x_j^\varepsilon}{t_j^\varepsilon}\right) e^{-i\frac{(x-x_j^\varepsilon)^2}{2\varepsilon t_j^\varepsilon}} + o(1). \tag{1.13}$$

*Moreover, we have $\limsup_{\varepsilon \to 0} \dfrac{t_j^\varepsilon}{\varepsilon} = +\infty$ for all $j \in \mathbb{N}$, and there is a $j \in \mathbb{N}$ such that $\limsup_{\varepsilon \to 0} t_j^\varepsilon \in [0,T]$.*

*Remark.* In the case when $x_j^\varepsilon$ and $t_j^\varepsilon$ have a limit as $\varepsilon$ goes to zero (for example if one has WKB initial data, in which case an immediate identification gives the value of those limits) then Theorem 1.4 gives the form of the solution $u^\varepsilon$ for all times.

*Remark.* From Theorem 1.1, (1.12) means that the obstruction for the initial data to be linearizable comes from $\Psi_\ell^\varepsilon$. Equation (1.13) asserts that this obstruction is due to quadratic oscillations. Notice that *a priori*, the profiles $\varphi_j$ do not belong to $H^1(\mathbb{R}^n)$. This information would require refined geometric assumptions on the initial data $u_0^\varepsilon$, an issue which we will not pursue here.

*Remark.* From a geometrical optics point of view, quadratic oscillations cause focusing at a point, which is the most degenerate caustic. In particular, if the initial data oscillates differently, a caustic may be formed, but the nonlinear term will remain negligible. This means that when the caustic is not a point, the scaling of the nonlinearity is sub-critical. This is very much in the spirit of the results of [JMR00], where this type of phenomenon is encountered, in a different setting.

Theorem 1.2 enables us to give a sufficient condition for condition (1.9) to be satisfied. Recall the definition of a Wigner measure. The Wigner measure of a family $(f^\varepsilon)_{0<\varepsilon\leq 1}$ bounded in $L^2(\mathbb{R}^n)$ is the weak limit (up to extraction) of its Wigner transform,

$$W^\varepsilon(f^\varepsilon)(x,\xi) = \int_{\mathbb{R}^n} f^\varepsilon\left(x - \frac{v\varepsilon}{2}\right) \overline{f^\varepsilon}\left(x + \frac{v\varepsilon}{2}\right) e^{i\xi \cdot v} \frac{dv}{(2\pi)^n}.$$

This limit is a positive Radon measure, whose study has proved to be efficient in semiclassical analysis, and in homogenization issues (see e.g. [LP93], [GMMP97]).

If $\mu$ and $\nu$ are two positive measures on the same measured space, the notation $\mu \perp \nu$ means that $\mu$ and $\nu$ are mutually singular, *i.e.* there exist measurable sets $A$ and $B$ such that $A \cap B = \emptyset$ and, for every measurable set $E$, $\mu(E) = \mu(E \cap A)$, $\nu(E) = \nu(E \cap B)$.

**Corollary 1.3.** *Let $T > 0$. Assume that for every Wigner measure $\mu_0$ associated with $u_0^\varepsilon$, for every $y \in \mathbb{R}^n$, for every $a \in [0, T]$, we have*

$$\mu_0(x, \xi) \perp \delta(x - y - a\xi) \otimes d\xi.$$

*Then (1.9) holds, and in particular, $v^\varepsilon$ is an approximation of $u^\varepsilon$ on the time interval $[0, T]$.*

We finally focus on the case where the initial data $u_0^\varepsilon$ have quadratic oscillations. We assume that

$$(1.14) \qquad u_0^\varepsilon(x) = \sum_{j=1}^{J} f_j(x) e^{-i\frac{|x-x_j|^2}{2\varepsilon t_j}} + r_0^\varepsilon(x),$$

where

$$(1.15) \qquad f_j \in \Sigma := \left\{ \phi \in H^1(\mathbb{R}^n);\ |x|\phi \in L^2(\mathbb{R}^n) \right\},$$

$x_j \in \mathbb{R}^n$, $t_j > 0$, with $(t_j, x_j) \neq (t_k, x_k)$ if $j \neq k$. We also assume that $r_0^\varepsilon$ is bounded in $H_\varepsilon^1$ and that its free evolution $r^\varepsilon$, defined by

$$(1.16) \qquad \begin{cases} i\varepsilon \partial_t r^\varepsilon + \frac{1}{2}\varepsilon^2 \Delta r^\varepsilon = 0, & (t,x) \in \mathbb{R}_+ \times \mathbb{R}^n, \\ r^\varepsilon_{|t=0} = r_0^\varepsilon, \end{cases}$$

satisfies (1.9) for some $T > 0$. For $1 \leq j \leq J$, we define $v_j^\varepsilon$ as the solution of the initial value problem

$$(1.17) \qquad \begin{cases} i\varepsilon \partial_t v_j^\varepsilon + \frac{1}{2}\varepsilon^2 \Delta v_j^\varepsilon = \varepsilon^{n\sigma} |v_j^\varepsilon|^{2\sigma} v_j^\varepsilon, & (t,x) \in \mathbb{R}_+ \times \mathbb{R}^n, \\ v_{j|t=0}^\varepsilon = f_j(x) e^{-i\frac{|x-x_j|^2}{2\varepsilon t_j}}. \end{cases}$$

Define $Z$ by

$$(1.18) \qquad Z = \mathcal{F} \circ S \circ \mathcal{F}^{-1},$$

where $S$ is the nonlinear scattering operator associated to (1.4) (see e.g. [Caz93]) and

$$(1.19) \qquad \mathcal{F}f(\xi) = \widehat{f}(\xi) = \frac{1}{(2i\pi)^{n/2}} \int_{\mathbb{R}^n} e^{-ix\cdot\xi} f(x) dx.$$

From [Car00], we know that for $1 \leq j \leq J$, the following asymptotics holds in $H_\varepsilon^1$ as $\varepsilon$ goes to zero,

$$(1.20) \qquad v_j^\varepsilon(t, x) = \begin{cases} \dfrac{1}{(1 - t/t_j)^{n/2}} f_j\left(\dfrac{x - x_j}{1 - t/t_j}\right) e^{i\frac{|x-x_j|^2}{2\varepsilon(t-t_j)}} + o(1), & \text{if } t < t_j, \\[2mm] \dfrac{1}{(t/t_j - 1)^{n/2}} Z f_j\left(\dfrac{x - x_j}{1 - t/t_j}\right) e^{i\frac{|x-x_j|^2}{2\varepsilon(t-t_j)}} + o(1), & \text{if } t > t_j. \end{cases}$$

Our final result is the following.

**Theorem 1.4.** *Assume that $u_0^\varepsilon$ is given by (1.14). Then for any $T > 0$ such that*

$$\limsup_{\varepsilon \to 0} \sup_{0 \leq t \leq T} \varepsilon^{n\sigma} \|r^\varepsilon(t)\|_{L^{2\sigma+2}}^{2\sigma+2} = 0,$$

*the following asymptotics holds in $L^\infty(0, T; L^2)$ as $\varepsilon$ goes to zero,*

$$u^\varepsilon = \sum_{j=1}^{J} v_j^\varepsilon + r^\varepsilon + o(1).$$

*Remark.* As we shall see in the proof of Theorem 1.4, one can prove a result in the space $L^\infty(0, T; H_\varepsilon^1)$ if one is prepared to suppose that $\varepsilon \nabla r^\varepsilon$ is linearizable. This is the case in particular if $r_0^\varepsilon \equiv 0$.

*Remark.* Compare for instance with the results of Oberguggenberger [Obe86], and Rauch and Reed [RR87]. They exhibited some cases where the solution of a nonlinear hyperbolic equation could be approximated as the sum of a solution to a nonlinear equation and of the free evolution of a linearizable data. As in our case, this was possible thanks to strong geometric assumptions. In nonlinear geometric optics, this phenomenon also occurs: in diffractive regime, an approximate solution is given by Schrödinger like equations. Lannes [Lan98] proved a splitting phenomenon in the propagation of different modes, also explained by a refined geometric analysis. This leads to a similar nonlinear superposition, whose statement and geometrical interpretation are similar to that of [BG99], [Ker01] or [Gal01].

The rest of the paper is organized as follows. In Section 2, we introduce the tools we will use throughout the paper, and in Section 3, we prove Theorem 1.1. In Section 4, we prove Theorem 1.2, and in Section 5, we prove Theorem 1.4. Finally in the Appendix, we consider the analog to Equation (1.6) with a harmonic potential, and we transpose the above results to that case.

*Acknowledgments.* The authors are indebted to Patrick Gérard for fruitful discussions on this work.

## 2. General results on NLS and applications

Before going into the proof of the theorems, let us recall some general tools and results on nonlinear Schrödinger equations, and their consequences for our study.

First, recall the classical definition (see e.g. [Caz93]),

**Definition 2.1.** A pair $(q, r)$ is **admissible** if $2 \leq r < \frac{2n}{n-2}$ (resp. $2 \leq r \leq \infty$ if $n = 1$, $2 \leq r < \infty$ if $n = 2$) and

$$\frac{2}{q} = \delta(r) := n\left(\frac{1}{2} - \frac{1}{r}\right).$$

Consider the initial value problem

(2.1) $$\begin{cases} i\partial_t \psi^\varepsilon + \frac{1}{2}\Delta \psi^\varepsilon = |\psi^\varepsilon|^{2\sigma} \psi^\varepsilon, & (t, x) \in \mathbb{R} \times \mathbb{R}^n, \\ \psi_{|t=0}^\varepsilon = \psi_0^\varepsilon(x). \end{cases}$$

Since $\sigma > 2/n$, and $\sigma < 2/(n-2)$ if $n \geq 3$, we know that if $\psi_0^\varepsilon \in H^1(\mathbb{R}^n)$, then (2.1) has a unique global solution $\psi^\varepsilon$, which satisfies

$$\psi^\varepsilon \in C(\mathbb{R}, H^1(\mathbb{R}^n)) \cap L^q(\mathbb{R}, W^{1,r}(\mathbb{R}^n)), \quad \forall (q, r) \text{ admissible}.$$

If $n \geq 3$, this result was proved by Ginibre and Velo ([GV85a], [GV85b], see also [Caz93]), and for $n = 1$ or $2$, it was proved by Nakanishi ([Nak99], see also [Nak01]).

Moreover, if $\psi_0^\varepsilon$ belong to a bounded domain of $H^1(\mathbb{R}^n)$, then $\psi^\varepsilon$ belongs to a bounded domain of $L^q(\mathbb{R}, W^{1,r}(\mathbb{R}^n))$ for any admissible pair $(q,r)$. Using the scaling

$$u^\varepsilon(t,x) = \frac{1}{\varepsilon^{n/2}} \psi^\varepsilon\left(\frac{t}{\varepsilon}, \frac{x}{\varepsilon}\right)$$

we get the following lemma.

**Lemma 2.2.** *The functions $u^\varepsilon$ and $v^\varepsilon$ satisfy the following properties.*

(1) $u^\varepsilon, v^\varepsilon \in C(\mathbb{R}; H^1(\mathbb{R}^n))$.

(2) *For any admissible pair $(q,r)$, there exists $C_r$ such that*

$$(2.2) \quad \|u^\varepsilon\|_{L^q(\mathbb{R};L^r)} + \|v^\varepsilon\|_{L^q(\mathbb{R};L^r)} + \|\varepsilon\nabla u^\varepsilon\|_{L^q(\mathbb{R};L^r)} + \|\varepsilon\nabla v^\varepsilon\|_{L^q(\mathbb{R};L^r)} \leq C_r \varepsilon^{-1/q}.$$

(3) *Recall that $U_0^\varepsilon(t) := e^{i\varepsilon \frac{t}{2}\Delta}$. For any admissible pair $(q,r)$, there exists $C_q$ independent of $\varepsilon$ such that*

$$\varepsilon^{\frac{1}{q}} \|U_0^\varepsilon(t)\varphi\|_{L^q(\mathbb{R};L^r)} \leq C_q \|\varphi\|_{L^2}.$$

(4) *For any admissible pairs $(q_1, r_1)$ and $(q_2, r_2)$ and any interval $I$, there exists $C_{r_1,r_2}$ independent of $\varepsilon$ and $I$ such that*

$$(2.3) \quad \varepsilon^{\frac{1}{q_1}+\frac{1}{q_2}} \left\|\int_{I\cap\{s\leq t\}} U_0^\varepsilon(t-s) F(s) ds \right\|_{L^{q_1}(I;L^{r_1})} \leq C_{r_1,r_2} \|F\|_{L^{q_2'}(I;L^{r_2'})}.$$

**Notation.** For $a^\varepsilon = a^\varepsilon(t,x)$ and $t > 0$, we write $\|a^\varepsilon\|_{L_t^q(L^r)} := \|a^\varepsilon\|_{L^q(0,t;L^r(\mathbb{R}^n))}$.

*Remark.* In the special case where $u_0^\varepsilon$ does not depend on $\varepsilon$, say $u_0^\varepsilon = f \in \mathcal{S}(\mathbb{R}^n)$ a function of the Schwartz space, it is easy to see, through stationary phase arguments, that

$$\|v^\varepsilon(t,.) - f\|_{L^2 \cap L^\infty} \xrightarrow[\varepsilon \to 0]{} 0, \; \forall t \geq 0.$$

In that case, the above estimate (2.2) is therefore far from being sharp, in terms of powers of $\varepsilon$. On the other hand, if we choose the same initial data as in [Car00], $u_0^\varepsilon(x) = f(x) e^{-ix^2/(2\varepsilon)}$, then the computations of [Car00] show that for $t > 1$, $\|v^\varepsilon\|_{L_t^q(L^r)}$ and $\|u^\varepsilon\|_{L_t^q(L^r)}$ are exactly of order $\varepsilon^{-1/q}$, when $(q,r)$ is admissible. This shows that Strichartz estimates are sharp in terms of powers of $\varepsilon$, and suggests that the best possible estimates are those of this case, which is that of the worst geometry from the viewpoint of geometric optics (all the rays meet at one point). Bearing this credo in mind can be helpful for the intuition in the proof of Theorem 1.1.

The well-known conservations of mass and energy for (2.1) yield,

**Lemma 2.3.** *The following quantities are independent of time.*

- *Mass:* $\|u^\varepsilon(t)\|_{L^2} = \|v^\varepsilon(t)\|_{L^2} = \|u_0^\varepsilon\|_{L^2}$.
- *Linear energy:*

$$E_0^\varepsilon(t) := \frac{1}{2} \|\varepsilon\nabla_x v^\varepsilon(t)\|_{L^2}^2 = E_0^\varepsilon(0).$$

- *Nonlinear energy:*

$$E^\varepsilon(t) := \frac{1}{2}\|\varepsilon\nabla_x u^\varepsilon(t)\|_{L^2}^2 + \frac{\varepsilon^{n\sigma}}{\sigma+1} \|u^\varepsilon(t)\|_{L^{2\sigma+2}}^{2\sigma+2} = E^\varepsilon(0).$$

In particular, since $u_0^\varepsilon$ is bounded in $H_\varepsilon^1$, $u^\varepsilon$ and $v^\varepsilon$ are bounded in $H_\varepsilon^1$, uniformly in time.

We conclude this section by stating the Gagliardo-Nirenberg inequalities we will use.

**Lemma 2.4.** *Let $r \geq 2$, and if $n \geq 3$, assume that $r < 2n/(n-2)$. Then there exists $C_r$ such that for any $f \in H^1(\mathbb{R}^n)$,*
$$\|f\|_{L^r} \leq C_r \|f\|_{L^2}^{1-\delta(r)} \|\nabla f\|_{L^2}^{\delta(r)}.$$
*In particular, for any $\varepsilon \in ]0,1]$,*
$$\|f\|_{L^r} \leq C_r \varepsilon^{-\delta(r)} \|f\|_{L^2}^{1-\delta(r)} \|\varepsilon \nabla f\|_{L^2}^{\delta(r)},$$
*and every solution bounded in $H_\varepsilon^1$ is bounded in $L_\varepsilon^{2\sigma+2}$.*

## 3. The linearizability condition

We can now prove Theorem 1.1. We first prove that (1) implies (2). Let
$$R := \limsup_{\varepsilon \to 0} \sup_{0 \leq t \leq T} \left| \frac{1}{2} \|\varepsilon \nabla_x u^\varepsilon(t)\|_{L^2}^2 + \frac{\varepsilon^{n\sigma}}{\sigma+1} \|u^\varepsilon(t)\|_{L^{2\sigma+2}}^{2\sigma+2} \right.$$
$$\left. - \frac{1}{2} \|\varepsilon \nabla_x v^\varepsilon(t)\|_{L^2}^2 - \frac{\varepsilon^{n\sigma}}{\sigma+1} \|v^\varepsilon(t)\|_{L^{2\sigma+2}}^{2\sigma+2} \right|.$$

On the one hand, assumption (1) implies that
$$R \leq \limsup_{\varepsilon \to 0} \sup_{0 \leq t \leq T} \frac{\varepsilon^{n\sigma}}{\sigma+1} \int \left| |u^\varepsilon(t,x)|^{2\sigma+2} - |v^\varepsilon(t,x)|^{2\sigma+2} \right| dx$$
$$\lesssim \limsup_{\varepsilon \to 0} \sup_{0 \leq t \leq T} \varepsilon^{n\sigma} \int \left( |u^\varepsilon(t,x)|^{2\sigma+1} + |v^\varepsilon(t,x)|^{2\sigma+1} \right) |u^\varepsilon(t,x) - v^\varepsilon(t,x)| dx.$$

Hölder's inequality yields
$$R \lesssim \limsup_{\varepsilon \to 0} \sup_{0 \leq t \leq T} \varepsilon^{n\sigma} \|u^\varepsilon(t) - v^\varepsilon(t)\|_{L^{2\sigma+2}} \left( \|u^\varepsilon(t)\|_{L^{2\sigma+2}} + \|v^\varepsilon(t)\|_{L^{2\sigma+2}} \right)^{2\sigma+1}.$$

From Lemma 2.3 and Lemma 2.4, this implies
$$R \lesssim \limsup_{\varepsilon \to 0} \sup_{0 \leq t \leq T} \varepsilon^{n\sigma} \|u^\varepsilon(t) - v^\varepsilon(t)\|_{L^{2\sigma+2}} \varepsilon^{-n\sigma \frac{2\sigma+1}{2\sigma+2}}$$
$$\lesssim \limsup_{\varepsilon \to 0} \sup_{0 \leq t \leq T} \varepsilon^{\delta(2\sigma+2)} \|u^\varepsilon(t) - v^\varepsilon(t)\|_{L^{2\sigma+2}}$$
$$\lesssim \limsup_{\varepsilon \to 0} \sup_{0 \leq t \leq T} \|u^\varepsilon(t) - v^\varepsilon(t)\|_{L^2}^{1-\delta(2\sigma+2)} \|\varepsilon \nabla u^\varepsilon(t) - \varepsilon \nabla v^\varepsilon(t)\|_{L^2}^{\delta(2\sigma+2)}.$$

We conclude from (1) that $R = 0$. On the other hand, the conservation of energy stated in Lemma 2.3 yield
$$R = \limsup_{\varepsilon \to 0} \sup_{0 \leq t \leq T} \frac{\varepsilon^{n\sigma}}{\sigma+1} \left| \|u_0^\varepsilon\|_{L^{2\sigma+2}}^{2\sigma+2} - \|v^\varepsilon(t)\|_{L^{2\sigma+2}}^{2\sigma+2} \right|,$$
and since $u_0^\varepsilon$ is supposed to go to zero in $L_\varepsilon^{2\sigma+2}$, this proves that (1) implies (2).

To prove that (2) implies (1), introduce the (expected) remainder $w^\varepsilon := u^\varepsilon - v^\varepsilon$. It solves,
$$(3.1) \quad \begin{cases} i\varepsilon \partial_t w^\varepsilon + \frac{1}{2} \varepsilon^2 \Delta w^\varepsilon = \varepsilon^{n\sigma} |u^\varepsilon|^{2\sigma} u^\varepsilon, & (t,x) \in \mathbb{R}_+ \times \mathbb{R}^n, \\ w^\varepsilon_{|t=0} = 0. \end{cases}$$

In a first step, we prove that
$$\sup_{0 \leq t \leq T} \|w^\varepsilon(t)\|_{L^2} \xrightarrow[\varepsilon \to 0]{} 0.$$

Writing $|u^\varepsilon|^{2\sigma} u^\varepsilon = |u^\varepsilon|^{2\sigma} w^\varepsilon + |u^\varepsilon|^{2\sigma} v^\varepsilon$ suggests the use of a Gronwall type argument, with $|u^\varepsilon|^{2\sigma} v^\varepsilon$ as a source term. Since in general (that is, when $n \geq 2$) we cannot expect $L^\infty$-estimates of $|u^\varepsilon|^{2\sigma}$, the linear term $|u^\varepsilon|^{2\sigma} w^\varepsilon$ must be handled with care. The following algebraic lemma, whose proof is left out, will help us.

**Lemma 3.1.** *Under the assumption $2/n < \sigma < 2/(n-2)$, one can find $\underline{q}, \underline{r}, \underline{s}$ and $\underline{k}$ such that*

$$\begin{cases} \dfrac{1}{\underline{r}'} = \dfrac{1}{\underline{r}} + \dfrac{2\sigma}{\underline{s}}, \\ \dfrac{1}{\underline{q}'} = \dfrac{1}{\underline{q}} + \dfrac{2\sigma}{\underline{k}}, \end{cases}$$

*and satisfying the additional conditions:*
- *The pair $(\underline{q}, \underline{r})$ is admissible,*
- $0 < \frac{2}{\underline{k}} < \delta(\underline{s}) < 1$.

Apply the Strichartz inequality (2.3) with $(q_1, r_1) = (\underline{q}, \underline{r})$, and $(q_2, r_2) = (\underline{q}, \underline{r})$ for the "linear" term, $(q_2, r_2) = (\infty, 2)$ for the source term. This yields,

$$\|w^\varepsilon\|_{L^{\underline{q}}_t(L^{\underline{r}})} \lesssim \varepsilon^{n\sigma - 1 - 2/\underline{q}} \||u^\varepsilon|^{2\sigma} w^\varepsilon\|_{L^{\underline{q}'}_t(L^{\underline{r}'})} + \varepsilon^{n\sigma - 1 - 1/\underline{q}} \||u^\varepsilon|^{2\sigma} v^\varepsilon\|_{L^1_t(L^2)}.$$

Hölder's inequality, along with Lemma 3.1, implies

(3.2) $\|w^\varepsilon\|_{L^{\underline{q}}_t(L^{\underline{r}})} \lesssim \varepsilon^{n\sigma - 1 - 2/\underline{q}} \|u^\varepsilon\|^{2\sigma}_{L^{\underline{k}}_t(L^{\underline{s}})} \|w^\varepsilon\|_{L^{\underline{q}}_t(L^{\underline{r}})} + \varepsilon^{n\sigma - 1 - 1/\underline{q}} \||u^\varepsilon|^{2\sigma} v^\varepsilon\|_{L^1_t(L^2)}.$

Our goal is to prove that up to increasing the constants (which we have not written so far), the left–hand side of (3.2) can be estimated by the second term of the right–hand side only.

**Lemma 3.2.** *Define $\theta_1$ by*

$$\theta_1 = \frac{\delta(\underline{s}) - 2/\underline{k}}{\delta(2\sigma + 2)}.$$

*Then $\theta_1 \in ]0, 1[$, and the pair $(q, r)$ defined as follows is admissible,*

$$\begin{cases} \dfrac{1}{\underline{s}} = \dfrac{1 - \theta_1}{r} + \dfrac{\theta_1}{2\sigma + 2}, \\ \dfrac{1}{\underline{k}} = \dfrac{1 - \theta_1}{q} + \dfrac{\theta_1}{\infty}. \end{cases}$$

Hölder's inequality then yields, along with Lemma 2.2,

$$\|u^\varepsilon\|^{2\sigma}_{L^{\underline{k}}_t(L^{\underline{s}})} \leq \|u^\varepsilon\|^{2\sigma(1-\theta_1)}_{L^q_t(L^r)} \|u^\varepsilon\|^{2\sigma\theta_1}_{L^\infty_t(L^{2\sigma+2})}$$
$$\lesssim \varepsilon^{-2\sigma(1-\theta_1)/q} \|u^\varepsilon\|^{2\sigma\theta_1}_{L^\infty_t(L^{2\sigma+2})}.$$

At this stage, we could estimate the last term thanks to Lemma 2.3; however, we could not absorb the first term of the right–hand side of (3.2) by the left–hand side, for the constant in factor of $\|w^\varepsilon\|_{L^{\underline{q}}_t(L^{\underline{r}})}$ need not be smaller than 1. We rather write

$$\|u^\varepsilon\|_{L^\infty_t(L^{2\sigma+2})} \leq \|v^\varepsilon\|_{L^\infty_t(L^{2\sigma+2})} + \|w^\varepsilon\|_{L^\infty_t(L^{2\sigma+2})}.$$

From (1.9),

$$\|v^\varepsilon\|_{L^\infty_t(L^{2\sigma+2})} = \varepsilon^{-\frac{n\sigma}{2\sigma+2}} o(1),$$

where $o(1)$ goes to zero as $\varepsilon$ goes to zero, uniformly for $t \in [0, T]$. From Lemma 2.4,

$$\|w^\varepsilon\|_{L^\infty_t(L^{2\sigma+2})} \lesssim \varepsilon^{-\frac{n\sigma}{2\sigma+2}} \|w^\varepsilon\|^{1-\delta(2\sigma+2)}_{L^\infty_t(L^2)} \|\varepsilon \nabla w^\varepsilon\|^{\delta(2\sigma+2)}_{L^\infty_t(L^2)}.$$

Gathering all these estimates, and using Lemma 2.2 to estimate $\varepsilon \nabla w^\varepsilon$, estimate (3.2) can be written as follows: there exist $C_*, C$ independent of $t \in [0, T]$ and $\varepsilon \in ]0, 1]$ such that, for any $t \in [0, T]$,

(3.3) $$\|w^\varepsilon\|_{L^{\underline{q}}_t(L^{\underline{r}})} \leq C_* \left( o(1) + \|w^\varepsilon\|^\gamma_{L^\infty_t(L^2)} \right) \|w^\varepsilon\|_{L^{\underline{q}}_t(L^{\underline{r}})}$$
$$+ C \varepsilon^{n\sigma - 1 - 1/\underline{q}} \||u^\varepsilon|^{2\sigma} v^\varepsilon\|_{L^1_t(L^2)},$$

where we denoted $\gamma = 2\sigma\theta_1(1-\delta(2\sigma+2))$. Notice that $\gamma > 0$; with our construction, it would not necessarily be so if we had $\sigma = 2/n$. Now, recall that we want to prove that $\|w^\varepsilon\|_{L^\infty_T(L^2)} = o(1)$. Thus for $\varepsilon$ sufficiently small, the term in factor of $\|w^\varepsilon\|_{L^q_t(L^r)}$ in the right–hand side of (3.3) should be less than, say, $1/2$, and the left–hand side could be estimated by the second term of the right–hand side. Let us make this argument rigorous. We know that at time $t = 0$, $w^\varepsilon = 0$. From Lemma 2.2, $w^\varepsilon$ is a continuous function of time with values in $L^2$, thus there exists $t_\varepsilon > 0$ such that
$$\|w^\varepsilon\|^\gamma_{L^\infty_{t_\varepsilon}(L^2)} < \frac{1}{4C_*}.$$

So long as

(3.4) $$\|w^\varepsilon\|^\gamma_{L^\infty_t(L^2)} < \frac{1}{4C_*}$$

holds, (3.3) thus yields
$$\|w^\varepsilon\|_{L^q_t(L^r)} \leq \left(C_* o(1) + \frac{1}{4}\right)\|w^\varepsilon\|_{L^q_t(L^r)} + C\varepsilon^{n\sigma - 1 - 1/q}\||u^\varepsilon|^{2\sigma}v^\varepsilon\|_{L^1_t(L^2)}.$$

Taking $\varepsilon$ sufficiently small, we have,
$$\|w^\varepsilon\|_{L^q_t(L^r)} \leq \frac{1}{2}\|w^\varepsilon\|_{L^q_t(L^r)} + C\varepsilon^{n\sigma - 1 - 1/q}\||u^\varepsilon|^{2\sigma}v^\varepsilon\|_{L^1_t(L^2)},$$

and we conclude that so long as (3.4) holds,

(3.5) $$\|w^\varepsilon\|_{L^q_t(L^r)} \lesssim \varepsilon^{n\sigma - 1 - 1/q}\||u^\varepsilon|^{2\sigma}v^\varepsilon\|_{L^1_t(L^2)}.$$

Applying Strichartz inequality (2.3) again, with now $(q_1, r_1) = (\infty, 2)$, yields, along with Lemma 3.1,

(3.6) $$\|w^\varepsilon\|_{L^\infty_t(L^2)} \lesssim \varepsilon^{n\sigma - 1 - 1/q}\|u^\varepsilon\|^{2\sigma}_{L^{\frac{k}{t}}_t(L^s)}\|w^\varepsilon\|_{L^q_t(L^r)} + \varepsilon^{n\sigma - 1}\||u^\varepsilon|^{2\sigma}v^\varepsilon\|_{L^1_t(L^2)}.$$

Now from Lemma 3.2, Lemma 2.2 and Lemma 2.3,
$$\varepsilon^{n\sigma - 1 - 1/q}\|u^\varepsilon\|^{2\sigma}_{L^{\frac{k}{t}}_t(L^s)} \lesssim \varepsilon^{1/q}.$$

Notice that at this stage, we do not decompose $u^\varepsilon$ as $v^\varepsilon + w^\varepsilon$. Using (3.5), we have: so long as (3.4) holds,

(3.7) $$\|w^\varepsilon\|_{L^\infty_t(L^2)} \lesssim \varepsilon^{n\sigma - 1}\||u^\varepsilon|^{2\sigma}v^\varepsilon\|_{L^1_t(L^2)}.$$

From Hölder's inequality,
$$\||u^\varepsilon|^{2\sigma}v^\varepsilon\|_{L^1_t(L^2)} \leq \|u^\varepsilon\|^{2\sigma}_{L^{2\sigma+1}_t(L^{4\sigma+2})}\|v^\varepsilon\|_{L^{2\sigma+1}_t(L^{4\sigma+2})}.$$

We use the analogue of Lemma 3.2.

**Lemma 3.3.** *Define $\theta_2$ by*
$$\theta_2 = \frac{(2\sigma + 2)(n\sigma - 2)}{n\sigma(2\sigma + 1)}.$$

*Then $\theta_2 \in ]0, 1[$, and the pair $(q_1, r_1)$ defined as follows is admissible,*
$$\begin{cases} \dfrac{1}{4\sigma + 2} = \dfrac{\theta_2}{2\sigma + 2} + \dfrac{1 - \theta_2}{r_1}, \\ \dfrac{1}{2\sigma + 1} = \dfrac{\theta_2}{\infty} + \dfrac{1 - \theta_2}{q_1}. \end{cases}$$

From Hölder's inequality, we now have
$$\|v^\varepsilon\|_{L_t^{2\sigma+1}(L^{4\sigma+2})} \leq \|v^\varepsilon\|_{L_t^\infty(L^{2\sigma+2})}^{\theta_2} \|v^\varepsilon\|_{L_t^{q_1}(L^{r_1})}^{1-\theta_2}.$$

Using (1.9) and Lemma 2.2, it follows that, since $\theta_2 > 0$,
$$\|v^\varepsilon\|_{L_t^{2\sigma+1}(L^{4\sigma+2})} \leq o(1)\varepsilon^{-\frac{n\sigma\theta_2}{2\sigma+2}}\varepsilon^{-(1-\theta_2)/q_1}.$$

For $u^\varepsilon$, we can replace $o(1)$ by a $O(1)$, from Lemma 2.3. Therefore, so long as (3.4) holds,
$$\|w^\varepsilon\|_{L_t^\infty(L^2)} \lesssim \varepsilon^{n\sigma-1} o(1) \varepsilon^{-(2\sigma+1)\left(\frac{n\sigma\theta_2}{2\sigma+2}+(1-\theta_2)/q_1\right)} = o(1).$$

Indeed, the powers of $\varepsilon$ cancel exactly, from Lemma 3.3. Taking $0 < \varepsilon \ll 1$, (3.4) thus holds up to time $t = T$, and the first part of (1) is proved:
$$\sup_{0 \leq t \leq T} \|w^\varepsilon(t)\|_{L^2} \xrightarrow[\varepsilon \to 0]{} 0.$$

To complete the proof of Theorem 1.1, we must show that the same holds for $\varepsilon \nabla w^\varepsilon$. We have all the tools we need to do so, and we will see that it is a consequence of the above asymptotics. Differentiating (3.1), we see that $\varepsilon \nabla w^\varepsilon$ solves the initial value problem,

(3.8) $$\begin{cases} i\varepsilon \partial_t (\varepsilon \nabla w^\varepsilon) + \frac{1}{2}\varepsilon^2 \Delta(\varepsilon \nabla w^\varepsilon) = \varepsilon^{n\sigma} \varepsilon \nabla \left(|u^\varepsilon|^{2\sigma} u^\varepsilon\right), & (t,x) \in \mathbb{R}_+ \times \mathbb{R}^n, \\ \varepsilon \nabla w^\varepsilon_{|t=0} = 0. \end{cases}$$

We can proceed with the same computations as before: notice that
$$\left|\varepsilon \nabla \left(|u^\varepsilon|^{2\sigma} u^\varepsilon\right)\right| \lesssim |u^\varepsilon|^{2\sigma} |\varepsilon \nabla w^\varepsilon| + |u^\varepsilon|^{2\sigma} |\varepsilon \nabla v^\varepsilon|.$$

The first term of the right–hand side plays the same role as $|u^\varepsilon|^{2\sigma} w^\varepsilon$ in the previous step. We can claim at this stage that the second term is small: the reason is that from the first step, we know that $u^\varepsilon$ goes to zero in $L_x^{2\sigma+2}$. Let us give some more convincing details.

Repeating the computations of the first step, but with a different admissible pair for the source term, $(\underline{q}, \underline{r})$ instead of $(\infty, 2)$, yields the analogue of (3.3),

(3.9) $$\begin{aligned} \|\varepsilon \nabla w^\varepsilon\|_{L_T^{\underline{q}}(L^{\underline{r}})} &\leq C'_* \left(o(1) + \|w^\varepsilon\|_{L_T^\infty(L^2)}^\gamma\right) \|\varepsilon \nabla w^\varepsilon\|_{L_T^{\underline{q}}(L^{\underline{r}})} \\ &\quad + C' \varepsilon^{n\sigma-1-2/\underline{q}} \left\||u^\varepsilon|^{2\sigma} \varepsilon \nabla v^\varepsilon\right\|_{L_T^{\underline{q}'}(L^{\underline{r}'})}, \end{aligned}$$

where we kept the notation $\gamma = 2\sigma\theta_1(1 - \delta(2\sigma+2))$, and possibly other constants $C'_*$ and $C'$. From the first step, if $\varepsilon$ is sufficiently small, then the factor $C'_*(o(1) + \|w^\varepsilon\|_{L_T^\infty(L^2)}^\gamma)$ is smaller than $1/2$, and (3.9) yields
$$\|\varepsilon \nabla w^\varepsilon\|_{L_T^{\underline{q}}(L^{\underline{r}})} \lesssim \varepsilon^{n\sigma-1-2/\underline{q}} \left\||u^\varepsilon|^{2\sigma} \varepsilon \nabla v^\varepsilon\right\|_{L_T^{\underline{q}'}(L^{\underline{r}'})}.$$

Plugging this estimate like in the first step, we have
$$\|\varepsilon \nabla w^\varepsilon\|_{L_T^\infty(L^2)} \lesssim \varepsilon^{n\sigma-1-1/\underline{q}} \left\||u^\varepsilon|^{2\sigma} \varepsilon \nabla v^\varepsilon\right\|_{L_T^{\underline{q}'}(L^{\underline{r}'})}.$$

From Lemma 3.1 and Hölder's inequality,
$$\left\||u^\varepsilon|^{2\sigma} \varepsilon \nabla v^\varepsilon\right\|_{L_T^{\underline{q}'}(L^{\underline{r}'})} \leq \|u^\varepsilon\|_{L_T^{\underline{k}}(L^{\underline{s}})}^{2\sigma} \|\varepsilon \nabla v^\varepsilon\|_{L_T^{\underline{q}}(L^{\underline{r}})}$$

Using Lemma 3.2 and Lemma 2.2, it follows
$$\begin{aligned} \left\||u^\varepsilon|^{2\sigma} \varepsilon \nabla v^\varepsilon\right\|_{L_T^{\underline{q}'}(L^{\underline{r}'})} &\lesssim \|u^\varepsilon\|_{L_T^q(L^r)}^{2\sigma(1-\theta_1)} \|u^\varepsilon\|_{L_T^\infty(L^{2\sigma+2})}^{2\sigma\theta_1} \varepsilon^{-1/\underline{q}} \\ &\lesssim \varepsilon^{-2\sigma(1-\theta_1)/q} \|u^\varepsilon\|_{L_T^\infty(L^{2\sigma+2})}^{2\sigma\theta_1} \varepsilon^{-1/\underline{q}}. \end{aligned}$$

But from (1.9) and the first step,

$$\|u^\varepsilon\|_{L_T^\infty(L^{2\sigma+2})} \leq \|v^\varepsilon\|_{L_T^\infty(L^{2\sigma+2})} + \|w^\varepsilon\|_{L_T^\infty(L^{2\sigma+2})}$$
$$\lesssim \varepsilon^{-\frac{n\sigma}{2\sigma+2}} o(1) + \varepsilon^{-\delta(2\sigma+2)} \|w^\varepsilon\|_{L_T^\infty(L^2)}^{1-\delta(2\sigma+2)} \|\varepsilon\nabla w^\varepsilon\|_{L_T^\infty(L^2)}^{\delta(2\sigma+2)}$$
$$\lesssim \varepsilon^{-\frac{n\sigma}{2\sigma+2}} o(1).$$

We conclude that $\|\varepsilon\nabla w^\varepsilon\|_{L_T^\infty(L^2)} = o(1)$, which completes the proof of Theorem 1.1.

*Remark.* Notice that in the proof of (1) $\Rightarrow$ (2), we did not use the assumption $\sigma > 2/n$, but only the conservations of mass and energy.

*Remark.* The proof of (2) $\Rightarrow$ (1) could have been achieved without knowing *a priori* the Strichartz estimates for $u^\varepsilon$, stated in Lemma 2.2, (2.2). Indeed, we exploited these estimates for only a finite number of $q$: since the knowledge of such estimates for $v^\varepsilon$ is straightforward, these estimates could have been obtained by a "so long" argument, as we used on the $L^\infty(L^2)$ norm of the remainder. Nevertheless, our proof relied at several steps on the assumption $\sigma > 2/n$, so knowing whether our result holds for instance when $\sigma = 2/n$ seems to be a challenging question.

## 4. When is the linearizability condition violated?

The goal of this section is to prove Theorem 1.2. We shall actually be proving a slightly more precise result than that stated in Theorem 1.2: it is in fact not necessary to assume that there is no focusing at time 0.

**Theorem 4.1.** *Let $T > 0$, suppose that $(H1), (H2)$ and $(H3)$ are satisfied, but assume that $(1.9)$ is not satisfied. Then up to the extraction of a subsequence, there exist two orthogonal families $(x_j^\varepsilon, t_j^\varepsilon)_{j\in\mathbb{N}}$ and $(\widetilde{x}_j^\varepsilon)_{j\in\mathbb{N}}$ in $\mathbb{R}^n \times \mathbb{R}_+$ and $\mathbb{R}^n$ respectively, a family $(\Psi_\ell^\varepsilon)_{\ell\in\mathbb{N}}$ bounded in $H_\varepsilon^1(\mathbb{R}^n)$, a family $(\Phi_j)_{j\in\mathbb{N}}$, bounded in $H^1(\mathbb{R}^n)$, and a family $(\varphi_j)_{j\in\mathbb{N}}$, bounded in $\mathcal{F}(H^1)(\mathbb{R}^n)$, such that:*

$$\text{(4.1)} \quad u_0^\varepsilon(x) = \Psi_\ell^\varepsilon(x) + \sum_{j=0}^\ell \frac{1}{\varepsilon^{\frac{n}{2}}} \Phi_j\left(\frac{x - \widetilde{x}_j^\varepsilon}{\varepsilon}\right) + r_\ell^\varepsilon(x),$$

$$\text{with } \limsup_{\varepsilon\to 0} \|U_0^\varepsilon(t) r_\ell^\varepsilon\|_{L^\infty(\mathbb{R}_+, L_\varepsilon^{2\sigma+2})} \xrightarrow[\ell\to\infty]{} 0,$$

*and for every $\ell \in \mathbb{N}$, the following asymptotics holds in $L^2(\mathbb{R}^n)$, as $\varepsilon \to 0$,*

$$\text{(4.2)} \quad \Psi_\ell^\varepsilon(x) = \sum_{j=0}^\ell \frac{1}{(t_j^\varepsilon)^{\frac{n}{2}}} \varphi_j\left(\frac{x - x_j^\varepsilon}{t_j^\varepsilon}\right) e^{-i\frac{(x-x_j^\varepsilon)^2}{2\varepsilon t_j^\varepsilon}} + o(1).$$

*We have also $\limsup_{\varepsilon\to 0} \dfrac{t_j^\varepsilon}{\varepsilon} = +\infty$.*

*In the case when $(H4)$ is satisfied, then for some $j$, $\limsup_{\varepsilon\to 0} t_j^\varepsilon \in [0, T]$; moreover $\Phi_j$ is equal to zero for all $j \in \mathbb{N}$, and the family $(\varphi_j)_{j\in\mathbb{N}}$ is nonempty.*

To prove Theorem 4.1, we shall mainly be using results of Keraani [Ker01], and in particular the profile decomposition of a sequence of solutions to the linear Schrödinger equation, bounded in energy. That result is only proved in dimension 3, but it clearly holds in all dimensions. In order to avoid unnecessary complications, and to use directly the theorem of [Ker01], we shall suppose in the proof of Theorem 4.1 that $n = 3$; once again that is purely for the sake of simplicity.

Let us start by rescaling our solution $v^\varepsilon$ in order to apply directly Theorem 1.6 of [Ker01]. We define the rescaled function

$$V^\varepsilon(s, y) := \varepsilon^{\frac{3}{2}} v^\varepsilon(\varepsilon s, \varepsilon y),$$

which satisfies the linear equation (1.8) with data

$$V_0^\varepsilon(y) := \varepsilon^{\frac{3}{2}} u_0^\varepsilon(\varepsilon y).$$

By $(H3)$, the family $V_0^\varepsilon$ is obviously bounded in $H^1$, so one can apply Theorem 1.6 of [Ker01]: we shall prove the following proposition.

**Proposition 4.2.** *Under the assumptions of Theorem 4.1, and up to the extraction of a subsequence, there exist two families $(s_j^\varepsilon)_{j\in\mathbb{N}}$ and $(y_j^\varepsilon)_{j\in\mathbb{N}}$ in $\mathbb{R}$ and $\mathbb{R}^3$ respectively, such that*

(4.3) $$\forall j \neq k, \quad \limsup_{\varepsilon\to 0} |s_j^\varepsilon - s_k^\varepsilon| + |y_j^\varepsilon - y_k^\varepsilon| = +\infty,$$

*and there exists a (non empty) bounded family $(V_j)_{j\in\mathbb{N}}$ of solutions to (1.8) in $H^1$, such that*

$$V^\varepsilon(s,y) = \sum_{j=0}^{\ell} V_j(s - s_j^\varepsilon, y - y_j^\varepsilon) + W_\ell^\varepsilon(s,y)$$

*with*

(4.4) $$\limsup_{\varepsilon\to 0} \|W_\ell^\varepsilon\|_{L^\infty(\mathbb{R}, L^{2\sigma+2})} = 0 \quad \text{when} \quad \ell \to \infty.$$

Let us prove that result. Applying directly Theorem 1.6 of [Ker01], we have the following decomposition for $V^\varepsilon$, up to an extraction:

(4.5) $$V^\varepsilon(s,y) = \sum_{j=0}^{\ell} \frac{1}{\sqrt{\eta_j^\varepsilon}} V_j\left(\frac{s - s_j^\varepsilon}{(\eta_j^\varepsilon)^2}, \frac{y - y_j^\varepsilon}{\eta_j^\varepsilon}\right) + W_\ell^\varepsilon(s,y)$$

where $\eta_j^\varepsilon \in \mathbb{R}_+ \setminus \{0\}$ are the scales of concentration, satisfying the following orthogonality condition:

$$\forall j \neq k, \quad \text{either} \quad \limsup_{\varepsilon\to 0} \frac{\eta_j^\varepsilon}{\eta_k^\varepsilon} + \frac{\eta_k^\varepsilon}{\eta_j^\varepsilon} = +\infty,$$

$$\text{or} \quad \eta_j^\varepsilon = \eta_k^\varepsilon \quad \text{and} \quad \limsup_{\varepsilon\to 0} \frac{|s_j^\varepsilon - s_k^\varepsilon| + |y_j^\varepsilon - y_k^\varepsilon|}{\eta_j^\varepsilon} = +\infty.$$

The remainder $W_\ell^\varepsilon$ satisfies

(4.6) $$\limsup_{\varepsilon\to 0} \|W_\ell^\varepsilon\|_{L^q(\mathbb{R}, L^r)} = 0 \quad \text{when} \quad \ell \to \infty,$$

and $\dfrac{2}{q} + \dfrac{3}{r} = \dfrac{1}{2}$, with $r < +\infty$. Such $(q,r)$ are said to be $\dot{H}^1$-admissible (as opposed to the $L^2$-admissible pairs of Definition 2.1, see [Ker01]).

Finally the $V_j$'s and $W_\ell^\varepsilon$ are solutions of (1.8) in $L^\infty(\mathbb{R}, \dot{H}^1)$.

So we need to prove the following facts.

- (1) Up to an extraction, $\eta_j^\varepsilon$ is equal to one for all $j$.
- (2) The $V_j$'s and $W_\ell^\varepsilon$ are bounded in $L^\infty(\mathbb{R}, H^1)$.
- (3) The limit (4.4) holds.
- (4) The family $(V_j)_{j\in\mathbb{N}}$ is not empty.

Note that the first result implies (2) and (3). Indeed suppose that $\eta_j^\varepsilon$ is equal to one for all $j$. Then by the orthogonality properties (4.3), one has

$$V^\varepsilon(s + s_j^\varepsilon, y + y_j^\varepsilon) \rightharpoonup V_j(s,y) \quad \text{in} \quad \mathcal{D}'(\mathbb{R} \times \mathbb{R}^3).$$

But we know by Lemma 2.3 that $V^\varepsilon$ is bounded in $L^\infty(\mathbb{R}, L^2)$, so it follows that for all $j \in \mathbb{N}$, the profiles $V_j$ are bounded in $L^\infty(\mathbb{R}, L^2)$; that implies also that

(4.7) $\quad (W_\ell^\varepsilon)_{0<\varepsilon\leq 1}$ is bounded in $L^\infty(\mathbb{R}, L^2)$, uniformly in $\ell \in \mathbb{N}$.

Since $L^2 \cap \dot H^1 = H^1$, point (2) follows. Then to find the limit (4.4) (point (3)), we simply use (4.6) with $q = +\infty$ and (4.7): the result follows by Hölder's inequality (using assumption $(H2)$).

Point (4) is due to the following observation: we have

$$\|V^\varepsilon\|_{L^\infty([0,T],L^{2\sigma+2})} \leq \sum_{j=1}^{\ell} \|V_j\|_{L^\infty([0,T],L^{2\sigma+2})} + \|W_\ell^\varepsilon\|_{L^\infty([0,T],L^{2\sigma+2})},$$

so since (1.9) is not satisfied, (4.4) implies that all of the $V_j$'s cannot be zero.

So to end the proof of Proposition 4.2, all we need to prove is that $\eta_j^\varepsilon$ is equal to one for all $j$. In fact we shall prove that, up to an extraction,

(4.8) $$\forall j \in \mathbb{N}, \quad \eta_j^\varepsilon = \lambda_j \in \mathbb{R}_+ \setminus \{0\},$$

and the expected result will follow simply by rescaling $V_j$ by $\lambda_j$. To prove (4.8), we are going to use the notion of $\varepsilon_n$–oscillating sequences. For more details on the subject, we refer to [BG99], [GMMP97].

**Definition 4.3.** *Let $(\varepsilon_n)_{n\in\mathbb{N}}$ be a given sequence in $\mathbb{R}_+ \setminus \{0\}$, and let $(V^n)$ be a bounded sequence in $\dot H^1$. Then the sequence $(V^n)$ is $\varepsilon_n$–oscillating if the following property holds:*

$$\limsup_{n\to\infty} \int_{\varepsilon_n|\xi|\leq R^{-1}} |\xi|^2 |\mathcal{F}(V^n)(\xi)|^2 \, d\xi + \int_{\varepsilon_n|\xi|\geq R} |\xi|^2 |\mathcal{F}(V^n)(\xi)|^2 \, d\xi \xrightarrow[R\to+\infty]{} 0.$$

*Remark.* For a time–dependent sequence $(V^n)$, uniformly bounded in $L^\infty(\mathbb{R}_+, \dot H^1)$, the definition holds taking the limit uniformly in time.

It is easy to see (see [Gér98] or [Ker01]) that $V^\varepsilon$ is $\eta_j^\varepsilon$–oscillating for every sequence $\eta_j^\varepsilon$ appearing in the decomposition (4.5). So the proof of (4.8) is a consequence of the following lemma.

**Lemma 4.4.** *Suppose the sequence $(V^\varepsilon)$ is $\eta^\varepsilon$–oscillating for some sequence $\eta^\varepsilon$. Then up to a subsequence, we can write $\eta^\varepsilon = \lambda$ for some $\lambda \in \mathbb{R}_+ \setminus \{0\}$.*

*Proof of Lemma 4.4.* Suppose that $(V^\varepsilon)$ is $\eta^\varepsilon$–oscillating. Then we can write, uniformly in time,

$$\|\nabla V^\varepsilon(s)\|_{L^2}^2 \lesssim \int_{R^{-1} \leq \eta^\varepsilon|\xi|\leq R} |\xi|^2 |\mathcal{F}(V^\varepsilon)|^2 \, d\xi + \delta(\varepsilon, R)$$

$$\lesssim \left(\frac{R}{\eta^\varepsilon}\right)^2 \|V^\varepsilon(s)\|_{L^2}^2 + \delta(\varepsilon, R),$$

where $\limsup_{\varepsilon\to 0} \delta(\varepsilon, R) \to 0$ as $R \to \infty$. The conservation of the energy yields

$$\|\nabla V^\varepsilon(s)\|_{L^2}^2 = \|\nabla V^\varepsilon(0)\|_{L^2}^2 = \|\varepsilon \nabla u_0^\varepsilon\|_{L^2}^2.$$

Up to a subsequence, we can suppose that this quantity is bounded from below by some $c > 0$ independent of $\varepsilon \in ]0,1]$ (otherwise, Condition (1.9) could not be violated, from Lemma 2.4). Fixing $R$ such that

$$\limsup_{\varepsilon\to 0} \delta(\varepsilon, R) \leq \frac{c}{2},$$

yields, up to an extraction,

$$\limsup_{\varepsilon\to 0} \eta^\varepsilon = \lambda \in \mathbb{R}_+.$$

Now suppose that $\lambda = 0$. One can write, for all time,

$$V^\varepsilon = V_R^\varepsilon + W_R^\varepsilon, \quad \text{with} \quad \mathcal{F}V_R^\varepsilon(t,\xi) := \mathbf{1}_{R^{-1}\leq \eta^\varepsilon|\xi|\leq R}\mathcal{F}V^\varepsilon(t,\xi),$$

and for all $\delta > 0$, if $R$ is large enough uniformly in $\varepsilon$ and $\eta^\varepsilon$, we have
$$\|W_R^\varepsilon\|_{L^\infty(\mathbb{R},H^1)} \leq \delta.$$
The Gagliardo–Nirenberg inequality (see Lemma 2.4) implies that $\|W_R^\varepsilon\|_{L^\infty(\mathbb{R},L^{2\sigma+2})}$ can be chosen arbitrarily small if $R$ is large enough, uniformly in $\varepsilon$ and $\eta^\varepsilon$. It follows that one can write

$$(4.9) \quad \begin{aligned} \|V^\varepsilon\|_{L^\infty([0,T],L^{2\sigma+2})}^{2\sigma+2} &\lesssim \|V_R^\varepsilon\|_{L^\infty([0,T],L^{2\sigma+2})}^{2\sigma+2} + o(1) \\ &\lesssim \|V_R^\varepsilon\|_{L^\infty([0,T],L^2)}^{(2\sigma+2)(1-\delta(2\sigma+2))} + o(1) \end{aligned}$$

where the second inequality is due again to Gagliardo–Nirenberg inequality and the boundedness of $V^\varepsilon$ in $H^1$.

Moreover, frequency localization once again implies that for all $s \in [0,T]$,
$$\|V_R^\varepsilon(s)\|_{L^2}^2 \lesssim (\eta^\varepsilon)^2 \int_{\frac{1}{R} \leq \eta^\varepsilon \xi \leq R} |\xi|^2 |\mathcal{F}V_R^\varepsilon(s,\xi)|^2 \, d\xi.$$

So the result follows, since by assumption the left–hand side in (4.9) does not go to zero and $(2\sigma + 2)(1 - \delta(2\sigma + 2)) = 2 - \sigma > 0$. Lemma 4.4 is proved (and with it, Proposition 4.2). □

Now let us finish the proof of Theorem 4.1. Using the decomposition given by Proposition 4.2 for $s = 0$, we get
$$V^\varepsilon(0,y) = \sum_{j=0}^{\ell} V_j(-s_j^\varepsilon, y - y_j^\varepsilon) + W_\ell^\varepsilon(0,y),$$
and going back to the definition of $v^\varepsilon$, it follows that one can write

$$(4.10) \quad u_0^\varepsilon(x) = \sum_{j=0}^{\ell} \frac{1}{\varepsilon^{\frac{3}{2}}} V_j\left(-\frac{t_j^\varepsilon}{\varepsilon}, \frac{x - x_j^\varepsilon}{\varepsilon}\right) + w_\ell^\varepsilon(x).$$

We have defined
$$(4.11) \quad t_j^\varepsilon := \varepsilon s_j^\varepsilon, \quad x_j^\varepsilon := \varepsilon y_j^\varepsilon \quad \text{and} \quad w_\ell^\varepsilon := \frac{1}{\varepsilon^{\frac{3}{2}}} W_\ell^\varepsilon\left(0, \frac{\cdot}{\varepsilon}\right).$$

Of course $U_0^\varepsilon(t)w_\ell^\varepsilon = \frac{1}{\varepsilon^{\frac{3}{2}}} W_\ell^\varepsilon\left(\frac{t}{\varepsilon}, \frac{\cdot}{\varepsilon}\right)$, so (4.4) can also be written

$$(4.12) \quad \limsup_{\varepsilon \to 0} \|U_0^\varepsilon(t)w_\ell^\varepsilon\|_{L^\infty(\mathbb{R},L_\varepsilon^{2\sigma+2})} \xrightarrow[\ell \to \infty]{} 0.$$

*Remark.* If one considers a family of initial data $r_0^\varepsilon$ going to zero in $L^2$ and uniformly bounded in $H_\varepsilon^1$, then clearly that convergence to zero holds uniformly in time for $U_0^\varepsilon(t)r_0^\varepsilon$ (due to Lemma 2.3), and Lemma 2.4 implies that $\|U_0^\varepsilon(t)r_0^\varepsilon\|_{L^\infty(\mathbb{R},L_\varepsilon^{2\sigma+2})}$ goes to zero.

Define the function
$$I_j^\varepsilon(t,x) := \frac{1}{\varepsilon^{\frac{3}{2}}} V_j\left(\frac{t - t_j^\varepsilon}{\varepsilon}, \frac{x - x_j^\varepsilon}{\varepsilon}\right).$$

We seek the asymptotic behavior of $I_j^\varepsilon(0,x)$ as $\varepsilon$ goes to zero. For every $j$, $V_j$ solves (1.8), we recall some classical results we will need on that equation (see e.g. [Rau91], [Caz93]).

**Lemma 4.5.** *Consider the initial value problem (1.8).*
• *If $V_0 \in L^2(\mathbb{R}^n)$, then (1.8) has a unique, global solution $V \in C(\mathbb{R}, L^2)$. Moreover, its asymptotic behavior as $t \to \pm\infty$ is given by*
$$\left\|V(t,x) - \frac{e^{i\frac{x^2}{2t}}}{t^{n/2}} \widehat{V_0}\left(\frac{x}{t}\right)\right\|_{L^2} \xrightarrow[t \to \pm\infty]{} 0,$$

where the Fourier transform is defined by (1.19), and we use the notation

$$\frac{1}{t^{n/2}} = \frac{i^n}{|t|^{n/2}} \quad \text{if } t < 0 .$$

- If $V_0 \in H^1(\mathbb{R}^n)$, then (1.8) has a unique, global solution $V \in C(\mathbb{R}, H^1)$. Moreover, the following quantities are independent of time,

$$\text{Mass: } \|V(t)\|_{L^2} = \|V_0\|_{L^2}$$
$$\text{Energy: } \|\nabla_x V(t)\|_{L^2} = \|\nabla V_0\|_{L^2}.$$

- If $V_0 \in \Sigma$, where $\Sigma$ is defined in (1.15), then (1.8) has a unique, global solution $V \in C(\mathbb{R}, \Sigma)$; it satisfies

$$\left\| V(t,x) - \frac{e^{i\frac{x^2}{2t}}}{t^{n/2}} \widehat{V_0}\left(\frac{x}{t}\right) \right\|_{H^1} \xrightarrow[t \to \pm\infty]{} 0.$$

For $\ell \in \mathbb{N}$, define

$$J_\ell := \left\{ j \in \{0, \ldots, \ell\}, \quad \limsup_{\varepsilon \to 0} \frac{t_j^\varepsilon}{\varepsilon} = +\infty \right\},$$
$$\widetilde{J}_\ell := \left\{ j \in \{0, \ldots, \ell\}, \quad \limsup_{\varepsilon \to 0} \frac{t_j^\varepsilon}{\varepsilon} \in \mathbb{R} \right\},$$

and we denote by $J_\infty$ (resp. $\widetilde{J}_\infty$) the union of all these $J_\ell$ (resp. $\widetilde{J}_\ell$).

The following proposition enables us to transform the profiles $V_j$ of (4.10) into the form given in Theorem 4.1.

**Proposition 4.6.** *Up to extraction of a subsequence, still denoted $\varepsilon$, we have:*
- *If $j \in \widetilde{J}_\infty$, then there exists $\Phi_j \in H^1$ such that, as $\varepsilon$ goes to zero, the following asymptotics holds in $H^1_\varepsilon$,*

$$I_j^\varepsilon(0,x) = \frac{1}{\varepsilon^{\frac{3}{2}}} \Phi_j\left(\frac{x - x_j^\varepsilon}{\varepsilon}\right) + o(1).$$

- *If $j \in \mathbb{N} \setminus \widetilde{J}_\infty$, then there exists $\psi_j \in \mathcal{F}(H^1)$ such that, as $\varepsilon$ goes to zero, the following asymptotics holds in $L^2$,*

$$I_j^\varepsilon(0,x) = \frac{1}{\left(t_j^\varepsilon\right)^{\frac{3}{2}}} e^{-\frac{i(x - x_j^\varepsilon)^2}{2\varepsilon t_j^\varepsilon}} \psi_j\left(\frac{x - x_j^\varepsilon}{t_j^\varepsilon}\right) + o(1).$$

- *For $j \in \mathbb{N} \setminus \widetilde{J}_\infty$, if we assume moreover that $V_j(0, \cdot) \in \Sigma$, then $\psi_j \in \Sigma$, and the above asymptotics holds not only in $L^2$, but also in $H^1_\varepsilon$.*
- *If $j \in \mathbb{N} \setminus (J_\infty \cup \widetilde{J}_\infty)$, then under the assumptions of Theorem 4.1,*

$$\limsup_{\varepsilon \to 0} \|I_j^\varepsilon\|_{L^\infty([0,T], L^{2\sigma+2}_\varepsilon)} = 0.$$

*Proof of Proposition 4.6.* By definition of $I_j^\varepsilon$, we have

$$I_j^\varepsilon(0,x) = \frac{1}{\varepsilon^{\frac{3}{2}}} V_j\left(-\frac{t_j^\varepsilon}{\varepsilon}, \frac{x - x_j^\varepsilon}{\varepsilon}\right).$$

If $j \in \widetilde{J}_\infty$, then the first point of Proposition 4.6 follows from the second point of Lemma 4.5, with

$$\Phi_j(x) = V_j(-\lambda_j, x),$$

where $\lambda_j := \limsup_{\varepsilon \to 0}(t_j^\varepsilon/\varepsilon)$.

If $j \in \mathbb{N} \setminus \widetilde{J}_\infty$, then the second point of the proposition follows from the first point of Lemma 4.5, with

$$\psi_j(x) = \widehat{V_j}(0, -x),$$

where the Fourier transform is taken with respect to the space variable only. Similarly, the third point of the proposition is a consequence of the last point of Lemma 4.5.

For the last point of Proposition 4.6, we need to consider more precisely the case when $\lambda_j = -\infty$. The idea of the result we want to prove is that if $\lambda_j = -\infty$, then for $\varepsilon$ small enough, we have $t_j^\varepsilon < 0$, so focusing for $I_j^\varepsilon$ has taken place at a negative time. Since there cannot be more than one focusing time, the solution is linearizable for positive times. Let us make that idea precise. For $\delta > 0$, there exists $F_j \in \Sigma$ such that
$$\|V_j(0, \cdot) - F_j\|_{H^1} \leq \delta.$$
Let $\widetilde{V}_j(t, x)$ be the solution of the linear Schrödinger equation (1.8) with initial data $F_j$. Then from the conservations of mass and energy, we have

(4.13) $$\|V_j - \widetilde{V}_j\|_{L^\infty(\mathbb{R}, H^1)} \leq \delta.$$

Define $\widetilde{I}_j^\varepsilon$ as the counterpart of $I_j^\varepsilon$ by
$$\widetilde{I}_j^\varepsilon(t, x) := \frac{1}{\varepsilon^{\frac{3}{2}}} \widetilde{V}_j\left(\frac{t - t_j^\varepsilon}{\varepsilon}, \frac{x - x_j^\varepsilon}{\varepsilon}\right).$$

From (4.13) and Lemma 2.4, we have,
$$\sup_{t \in [0, T]} \|I_j^\varepsilon(t, \cdot) - \widetilde{I}_j^\varepsilon(t, \cdot)\|_{L_\varepsilon^{2\sigma+2}} \leq C_{2\sigma+2} \delta.$$

The uniformity in time stems from the fact that for every $t \in [0, T]$,
$$\frac{|t_j^\varepsilon - t|}{\varepsilon} = \frac{t - t_j^\varepsilon}{\varepsilon} \geq \frac{-t_j^\varepsilon}{\varepsilon} \xrightarrow[\varepsilon \to 0]{} +\infty,$$
since we assumed $\lambda_j = -\infty$.

On the other hand, the last part of Lemma 4.5 implies $\widetilde{V}_j \in C(\mathbb{R}, \Sigma)$, and
$$\widetilde{I}_j^\varepsilon(t, x) = \frac{1}{(t - t_j^\varepsilon)^{\frac{3}{2}}} e^{i \frac{|x - x_j^\varepsilon|^2}{2\varepsilon(t - t_j^\varepsilon)}} \widehat{F}_j\left(\frac{x - x_j^\varepsilon}{t - t_j^\varepsilon}\right) + o(1), \text{ in } L^\infty([0, T], H_\varepsilon^1).$$

This explicit expression and Lemma 2.4 yield

(4.14) $$\|\widetilde{I}_j^\varepsilon\|_{L^\infty([0,T], L_\varepsilon^{2\sigma+2})} \leq \sup_{t \in [0,T]} \left|\frac{\varepsilon}{t_j^\varepsilon - t}\right|^{\frac{3\sigma}{2\sigma+2}} \|\widetilde{F}_j\|_{L^{2\sigma+2}} + o(1).$$

We assumed $\lambda_j = -\infty$, thus for $\varepsilon$ small enough and all $t \in [0, T]$,
$$\frac{\varepsilon}{|t_j^\varepsilon - t|} = \frac{\varepsilon}{t - t_j^\varepsilon} \leq \frac{\varepsilon}{-t_j^\varepsilon} \xrightarrow[\varepsilon \to 0]{} 0.$$

Therefore,
$$\limsup_{\varepsilon \to 0} \|I_j^\varepsilon\|_{L^\infty([0,T], L_\varepsilon^{2\sigma+2})} \leq C_{2\sigma+2} \delta,$$
and since $\delta > 0$ is arbitrary, this completes the proof of Proposition 4.6. $\square$

Now let us deduce the theorem. Write
$$u_0^\varepsilon(x) = \sum_{j \in J_\ell} I_j^\varepsilon(0, x) + \sum_{j \in \widetilde{J}_\ell} I_j^\varepsilon(0, x) + \rho_\ell^\varepsilon(x),$$

The remainder $\rho_\ell^\varepsilon$ satisfies the limit (4.1): it is indeed the sum of $w_\ell^\varepsilon$, which satisfies (4.12), and of the profiles associated with $\lambda_j = -\infty$, which satisfy the desired limit by the last point of Proposition 4.6.

For $j \in \widetilde{J}_\infty$, define $\widetilde{x}_j^\varepsilon = x_j^\varepsilon$. From the first point of Proposition 4.6, we also know that there exists a family $(\Phi_j)$ in $H^1$ such that for any $\ell \in \mathbb{N}$,

$$(4.15) \quad u_0^\varepsilon(x) = \sum_{j \in J_\ell} I_j^\varepsilon(0, x) + \sum_{j \in \widetilde{I}_\ell} \frac{1}{\varepsilon^{\frac{3}{2}}} \Phi_j \left( \frac{x - \widetilde{x}_j^\varepsilon}{\varepsilon} \right) + r_\ell^\varepsilon(x),$$

where $r_\ell^\varepsilon$ satisfies the limit (4.1), as the sum of $\rho_\ell^\varepsilon$ and of a term which is small in $H_\varepsilon^1$. From the second point of Proposition 4.6, we also know that there exists a family $(\varphi_j)$ in $L^2$ such that for any $\ell \in \mathbb{N}$, the following asymptotics holds in $L^2$ as $\varepsilon \to 0$,

$$\sum_{j \in J_\ell} I_j^\varepsilon(0, x) = \sum_{j \in J_\ell} \frac{1}{(t_j^\varepsilon)^{\frac{3}{2}}} e^{-\frac{i(x - x_j^\varepsilon)^2}{2\varepsilon t_j^\varepsilon}} \varphi_j \left( \frac{x - x_j^\varepsilon}{t_j^\varepsilon} \right) + o(1).$$

The orthogonality property on $(t_j^\varepsilon, x_j^\varepsilon)$ and on $(\widetilde{x}_j^\varepsilon)$ is due to definitions (4.11) and (4.3). The family $(\varphi_j)_{j \in \mathbb{N}}$ is the family $(\psi_j)_{j \in \mathbb{N}}$ from which we have removed the profiles associated to $\lambda_j = -\infty$. We used the obvious convention that $\psi_j \equiv 0$ if $j \in \widetilde{J}_\infty$, and similarly, we set $\Phi_j \equiv 0$ if $j \in \mathbb{N} \setminus \widetilde{J}_\infty$.

Now to end the proof of Theorem 4.1, we are left with the following proposition.

**Proposition 4.7.** *Under the assumptions of Th. 4.1, if additionally $(H4)$ is satisfied, then in decomposition (4.15) we have $\Phi_j = 0$ for every $j$, and the family $(\varphi_j)_{j \in \mathbb{N}}$ is nonempty. Moreover, there is an integer $j \in J_\infty$ such that $\limsup\limits_{\varepsilon \to 0} t_j^\varepsilon \in [0, T]$.*

Let us prove the result. It relies on the following lemma, which itself uses the orthogonality of the $\widetilde{x}_j^\varepsilon$'s. We shall postpone its proof to the end of this section.

**Lemma 4.8.** *Fix $\ell \in \mathbb{N}$. If for $1 \leq j \leq \ell$, $\Phi_j \in H^1$, then we have,*

$$\left\| \sum_{j=0}^\ell \frac{1}{\varepsilon^{\frac{3}{2}}} \Phi_j \left( \frac{\cdot - \widetilde{x}_j^\varepsilon}{\varepsilon} \right) \right\|_{L_\varepsilon^{2\sigma+2}}^{2\sigma+2} \xrightarrow[\varepsilon \to 0]{} \sum_{j=0}^\ell \|\Phi_j\|_{L^{2\sigma+2}}^{2\sigma+2}.$$

Now let us prove the first part of Proposition 4.7. Let $\delta > 0$. For every $j \in J_\infty$, there exists $F_j \in \Sigma$ such that

$$\|V_j(0, \cdot) - F_j\|_{H^1} \leq \frac{\delta}{2^{j+1}}.$$

Let $\widetilde{V}_j(t, x)$ be the solution of the linear Schrödinger equation (1.8) with initial data $F_j$. Then from the conservations of mass and energy, we have

$$(4.16) \quad \|V_j - \widetilde{V}_j\|_{L^\infty(\mathbb{R}, H^1)} \leq \frac{\delta}{2^{j+1}}.$$

Define $\widetilde{I}_j^\varepsilon$ as the counterpart of $I_j^\varepsilon$ by

$$\widetilde{I}_j^\varepsilon(t, x) := \frac{1}{\varepsilon^{\frac{3}{2}}} \widetilde{V}_j \left( \frac{t - t_j^\varepsilon}{\varepsilon}, \frac{x - x_j^\varepsilon}{\varepsilon} \right).$$

From (4.16) and Lemma 2.4, we have, for every $j \in J_\infty$,

$$\|I_j^\varepsilon(0, \cdot) - \widetilde{I}_j^\varepsilon(0, \cdot)\|_{L_\varepsilon^{2\sigma+2}} \leq C_{2\sigma+2} \frac{\delta}{2^{j+1}}.$$

Fix $\ell \in \mathbb{N}$. From (4.15), we have,

$$(4.17) \quad \left\| \sum_{j \in \widetilde{J}_\ell} I_j^\varepsilon(0, \cdot) \right\|_{L_\varepsilon^{2\sigma+2}} \leq \sum_{j \in J_\ell} \|\widetilde{I}_j^\varepsilon(0, \cdot)\|_{L_\varepsilon^{2\sigma+2}} + \|u_0^\varepsilon\|_{L_\varepsilon^{2\sigma+2}} + \|r_0^\varepsilon\|_{L_\varepsilon^{2\sigma+2}} + C_{2\sigma+2}\delta.$$

Proposition 4.6 yields profiles $\Phi_j \in H^1$ associated to $V_j$ for $j \in \widetilde{J}_\infty$, and $\widetilde{\varphi}_j \in \Sigma$ associated to $\widetilde{V}_j$ for $j \in J_\infty$. The asymptotics in $H^1_\varepsilon$ imply in particular, from Lemma 2.4,

$$\sum_{j \in J_\ell} \left\| \widetilde{I}^\varepsilon_j(0, \cdot) \right\|_{L^{2\sigma+2}_\varepsilon} = \sum_{j \in J_\ell} \left\| \frac{1}{(t^\varepsilon_j)^{\frac{3}{2}}} e^{-\frac{i(\cdot - x^\varepsilon_j)^2}{2\varepsilon t^\varepsilon_j}} \widetilde{\varphi}_j \left( \frac{\cdot - x^\varepsilon_j}{t^\varepsilon_j} \right) \right\|_{L^{2\sigma+2}_\varepsilon} + o(1), \text{ as } \varepsilon \to 0.$$

But we have

$$\left\| \frac{1}{(t^\varepsilon_j)^{\frac{3}{2}}} e^{-\frac{i(\cdot - x^\varepsilon_j)^2}{2\varepsilon t^\varepsilon_j}} \widetilde{\varphi}_j \left( \frac{\cdot - x^\varepsilon_j}{t^\varepsilon_j} \right) \right\|^{2\sigma+2}_{L^{2\sigma+2}_\varepsilon} = \left( \frac{\varepsilon}{t^\varepsilon_j} \right)^{3\sigma} \|\widetilde{\varphi}_j\|^{2\sigma+2}_{L^{2\sigma+2}},$$

so by definition of $J_\ell$ we get

$$\limsup_{\varepsilon \to 0} \left\| \frac{1}{(t^\varepsilon_j)^{\frac{3}{2}}} e^{-\frac{i(\cdot - x^\varepsilon_j)^2}{2\varepsilon t^\varepsilon_j}} \widetilde{\varphi}_j \left( \frac{\cdot - x^\varepsilon_j}{t^\varepsilon_j} \right) \right\|_{L^{2\sigma+2}_\varepsilon} = 0.$$

Then (4.1), (H4) and (4.17) imply that

$$\limsup_{\ell \to \infty} \limsup_{\varepsilon \to 0} \left\| \sum_{j \in \widetilde{J}_\ell} \frac{1}{\varepsilon^{\frac{3}{2}}} \Phi_j \left( \frac{\cdot - \widetilde{x}^\varepsilon_j}{\varepsilon} \right) \right\|_{L^{2\sigma+2}_\varepsilon} \leq C_{2\sigma+2} \delta.$$

In particular, for $\ell \geq \ell_0$ sufficiently large,

(4.18) $$\limsup_{\varepsilon \to 0} \left\| \sum_{j \in \widetilde{J}_\ell} \frac{1}{\varepsilon^{\frac{3}{2}}} \Phi_j \left( \frac{\cdot - \widetilde{x}^\varepsilon_j}{\varepsilon} \right) \right\|_{L^{2\sigma+2}_\varepsilon} \leq 2 C_{2\sigma+2} \delta.$$

Lemma 4.8 along with the above estimate imply that for every $\ell \in \mathbb{N}$ and every $j \in \widetilde{I}_\ell$,

$$\|\Phi_j\|_{L^{2\sigma+2}} \leq 2 C_{2\sigma+2} \delta.$$

Since $\delta > 0$ is arbitrary, this means that all the $\Phi_j$'s are zero. But Proposition 4.2 states that the family $V_j$ is non empty, so the $\varphi_j$'s cannot all be equal to zero. That yields the first part of Proposition 4.7.

Let us now prove the second part. We apply the operator $U^\varepsilon_0(t)$ to decomposition (4.15), and using the fact that all the $\Phi_j$'s are zero, we get for $\ell$ large enough, for all times $t \in [0, T]$,

$$\|U^\varepsilon_0(t) u^\varepsilon_0\|_{L^\infty([0,T], L^{2\sigma+2}_\varepsilon)} \leq \sum_{j \in I_\ell} \|I^\varepsilon_j\|_{L^\infty([0,T], L^{2\sigma+2}_\varepsilon)} + \delta(T, \ell, \varepsilon),$$

where $\limsup_{\varepsilon \to 0} \delta(T, \ell, \varepsilon) \to 0$ as $\ell \to \infty$. Now suppose that for all integers $j \in J_\infty$,

(4.19) $$\limsup_{\varepsilon \to 0} t^\varepsilon_j > T.$$

This contradicts the conclusion of the proposition. Reasoning as in the proof of Proposition 4.6 (see in particular (4.14)), and using the inequality, for $t \in [0, T]$,

$$\frac{\varepsilon}{|t^\varepsilon_j - t|} = \frac{\varepsilon}{t^\varepsilon_j - t} \leq \frac{\varepsilon}{t^\varepsilon_j - T},$$

we get, for every $j \in J_\infty$,

$$\limsup_{\varepsilon \to 0} \|I^\varepsilon_j\|_{L^\infty([0,T], L^{2\sigma+2}_\varepsilon)} = 0.$$

The contradiction follows, since the principal assumption of Theorem 4.1 is that (1.9) is not satisfied. So finally the theorem is proved, up to the proof of Lemma 4.8.

*Proof of Lemma 4.8.* This lemma is very classical, and uses the orthogonality of the $\widetilde{x}_j^\varepsilon$'s (see for instance [Gér98]). The only difference with the usual case is that the nonlinearity is not polynomial; however it is $C^1$ (see also Section 5). So let $\ell$ be fixed, and let $j \neq k$ be in $\{1, .., \ell\}$. An immediate induction shows that is enough to prove that for all $j \neq k$,

$$\limsup_{\varepsilon \to 0} \varepsilon^{3\sigma} \int \frac{1}{\varepsilon^{\frac{3}{2}}} \left|\Phi_j\left(\frac{x - \widetilde{x}_j^\varepsilon}{\varepsilon}\right)\right| \frac{1}{\varepsilon^{\frac{3}{2}(2\sigma+1)}} \left|\Phi_k\left(\frac{x - \widetilde{x}_k^\varepsilon}{\varepsilon}\right)\right|^{2\sigma+1} dx = 0.$$

But a change of variables yields

$$\varepsilon^{3\sigma} \int \frac{1}{\varepsilon^{\frac{3}{2}}} \left|\Phi_j\left(\frac{x - \widetilde{x}_j^\varepsilon}{\varepsilon}\right)\right| \frac{1}{\varepsilon^{\frac{3}{2}(2\sigma+1)}} \left|\Phi_k\left(\frac{x - \widetilde{x}_k^\varepsilon}{\varepsilon}\right)\right|^{2\sigma+1} dx$$
$$= \int |\Phi_j|\left(y + \frac{\widetilde{x}_k^\varepsilon - \widetilde{x}_j^\varepsilon}{\varepsilon}\right) |\Phi_k|^{2\sigma+1}(y)\, dy,$$

and the result simply follows from the orthogonality of the family $(\widetilde{x}_j^\varepsilon)$, since at this stage, one can suppose that $\Phi_j$ is compactly supported. □

That concludes the proof of Theorem 4.1, hence of Theorem 1.2.

## 5. INITIAL DATA WITH QUADRATIC OSCILLATIONS

In this section, we aim at proving Theorem 1.4. We set

$$v^\varepsilon(t, x) := \sum_{j=1}^J v_j^\varepsilon(t, x).$$

We first recall an estimate obtained in [Car00].

**Lemma 5.1.** *Let $r \geq 2$, with $r < 2n/(n-2)$ if $n \geq 3$. For any $j \in [1, J]$, there exists $C = C(j, r)$ such that for any $t \geq 0$,*

$$(5.1) \qquad \|v_j^\varepsilon(t)\|_{L^r} \leq C \left(\frac{1}{|t - t_j| + \varepsilon}\right)^{\delta(r)}.$$

*Proof of Lemma 5.1.* This estimate follows from the results of [Car00], and we briefly recall its proof. First, the conservation of energy (stated in Lemma 2.3) shows that $v_j^\varepsilon$ is bounded in $H_\varepsilon^1$, and from Lemma 2.4, there exists $C = C(j, r)$ such that

$$\|v_j^\varepsilon(t)\|_{L^r} \leq C\varepsilon^{-\delta(r)}.$$

Introduce the Galilean operator $J_j^\varepsilon(t) = (x - x_j)/\varepsilon + i(t - t_j)\nabla_x$. Then from the pseudo-conformal conservation law (see e.g. [Caz93]), $J_j^\varepsilon(t)v_j^\varepsilon$ remains bounded in $L^2$ (this is straightforward in our case $\sigma > 2/n$). On the other hand, $J_j^\varepsilon$ also writes

$$J_j^\varepsilon(t) = i(t - t_j) e^{i\frac{|x - x_j|^2}{2\varepsilon(t - t_j)}} \nabla_x \left(e^{-i\frac{|x - x_j|^2}{2\varepsilon(t - t_j)}} \cdot\right),$$

so Lemma 2.4 yields

$$\|v_j^\varepsilon(t)\|_{L^r} \leq C|t - t_j|^{-\delta(r)},$$

which completes the proof of Lemma 5.1. □

We want to prove that the remainder defined by $w^\varepsilon = u^\varepsilon - r^\varepsilon - v^\varepsilon$ remains small in $L^\infty(0, T, L^2)$. Notice that $w^\varepsilon$ satisfies the following Schrödinger equation:

$$i\varepsilon \partial_t w^\varepsilon + \frac{\varepsilon^2}{2}\Delta w^\varepsilon = \varepsilon^{n\sigma}(f^\varepsilon + s^\varepsilon),$$

with

$$f^\varepsilon = |u^\varepsilon|^{2\sigma} u^\varepsilon - |v^\varepsilon + r^\varepsilon|^{2\sigma}(v^\varepsilon + r^\varepsilon),$$

$$s^\varepsilon = |v^\varepsilon + r^\varepsilon|^{2\sigma}(v^\varepsilon + r^\varepsilon) - \sum_{j=1}^{J} |v_j^\varepsilon|^{2\sigma} v_j^\varepsilon,$$

and is zero at time $t = 0$. The key estimate is the Strichartz inequality (2.3), which gives on any bounded interval $I$

$$(5.2) \quad \|w^\varepsilon\|_{L^{q_1}(I;L^{r_1})} \lesssim \varepsilon^{n\sigma-1-\frac{1}{q_1}-\frac{1}{q_2}} \|f^\varepsilon\|_{L^{q_2'}(I,L^{r_2'})} + \varepsilon^{n\sigma-1-\frac{1}{q_1}-\frac{1}{q_3}} \|s^\varepsilon\|_{L^{q_3'}(I,L^{r_3'})},$$

for all admissible pairs $(q_1, r_1)$, $(q_2, r_2)$, $(q_3, r_3)$. Since

$$|f^\varepsilon| \lesssim \left(|v^\varepsilon|^{2\sigma} + |u^\varepsilon|^{2\sigma} + |r^\varepsilon|^{2\sigma}\right) |w^\varepsilon|,$$

this term will be treated as a Gronwall type term, while $s^\varepsilon$ is a source term. This resumes the spirit of the proof of Th. 1.1.

We will proceed as follows, decomposing $[0, T]$ as the union of different types of intervals $I$ and arguing differently with regards to $I$. First, we examine the Gronwall type term, $\|f^\varepsilon\|_{L^{q_2'}(I,L^{r_2'})}$. Then, we study the source term, $\|s^\varepsilon\|_{L^{q_3'}(I,L^{r_3'})}$. Finally, we can conclude on $\|w^\varepsilon\|_{L^\infty(I;L^2)}$.

**The Gronwall type term**: For $\Lambda > 0$, we denote by $I_\Lambda$ some interval included in $\{|t - t_j| > \Lambda\varepsilon, \; \forall j \in \{1, ..., J\}\}$, and by $I_\eta$ for $\eta > 0$ any interval of length $\eta\varepsilon$ disjoint from all the intervals $I_\Lambda$.

**Lemma 5.2.** Let $(\underline{q}, \underline{r})$ be as in Lemma 3.1.
1) Consider $\delta > 0$. There exist $\eta_0 > 0$ and $\varepsilon_0 > 0$, such that for all $\eta \leq \eta_0$ and $0 < \varepsilon \leq \varepsilon_0$,

$$\|f^\varepsilon\|_{L^{\underline{q}'}(I_\eta, L^{\underline{r}'})} \leq \delta \varepsilon^{1-n\sigma+2/\underline{q}} \|w^\varepsilon\|_{L^{\underline{q}}(I_\eta, L^{\underline{r}})}.$$

2) There exist $C, \gamma > 0$ such that for any $\delta > 0$, there exist $\Lambda_0 \geq 1$ and $\varepsilon_1 > 0$, such that for all $\Lambda \geq \Lambda_0$ and $0 < \varepsilon \leq \varepsilon_1$,

$$\|f^\varepsilon\|_{L^{\underline{q}'}(I_\Lambda, L^{\underline{r}'})} \leq C\varepsilon^{1-n\sigma+2/\underline{q}} \left(\delta + \|w^\varepsilon\|_{L^\infty(I_\Lambda, L^2)}^\gamma\right) \|w^\varepsilon\|_{L^{\underline{q}}(I_\Lambda, L^{\underline{r}})}.$$

*Proof of Lemma 5.2.* Using Hölder's inequality and Lemma 3.1, we obtain for any interval $I$,

$$(5.3) \quad \|f^\varepsilon\|_{L^{\underline{q}'}(I,L^{\underline{r}'})} \lesssim \left(\|u^\varepsilon\|_{L^{\underline{k}}(I,L^{\underline{s}})}^{2\sigma} + \|v^\varepsilon\|_{L^{\underline{k}}(I,L^{\underline{s}})}^{2\sigma} + \|r^\varepsilon\|_{L^{\underline{k}}(I,L^{\underline{s}})}^{2\sigma}\right) \|w^\varepsilon\|_{L^{\underline{q}}(I,L^{\underline{r}})}.$$

**1)** Let us study the case of an interval of the form $I_\eta$ for some $\eta > 0$. The Gagliardo-Nirenberg inequality of Lemma 2.4 yields for any $t \in I_\eta$,

$$\|u^\varepsilon(t)\|_{L^{\underline{s}}} \lesssim \varepsilon^{-\delta(\underline{s})} \|u^\varepsilon(t)\|_{L^2}^{1-\delta(\underline{s})} \|\varepsilon\nabla_x u^\varepsilon(t)\|_{L^2}^{\delta(\underline{s})},$$

and similarly for $v^\varepsilon$ and $r^\varepsilon$. Therefore, using the boundedness of the $H_\varepsilon^1$-norm of the solutions of Schrödinger equations, we obtain

$$\|u^\varepsilon\|_{L^{\underline{k}}(I_\eta, L^{\underline{s}})}^{2\sigma} \lesssim \left(\int_{I_\eta} \varepsilon^{-\delta(\underline{s})\underline{k}} dt\right)^{\frac{2\sigma}{\underline{k}}}$$

$$\lesssim \eta^{2\sigma/\underline{k}} \varepsilon^{2\sigma(\frac{1}{\underline{k}} - \delta(\underline{s}))}.$$

Since $2\sigma\left(\delta(\underline{s}) - \frac{1}{\underline{k}}\right) = -\frac{2}{\underline{q}} + n\sigma - 1$, we obtain **1)** in Lemma 5.2.

**2)** Let us consider now intervals of the type $I_\Lambda$ for some $\Lambda > 0$. Using Lemma 3.2, Lemma 2.2 and Hölder's inequality, we obtain

$$\|u^\varepsilon\|^{2\sigma}_{L^{\underline{k}}(I_\Lambda, L^{\underline{s}})} \lesssim \|u^\varepsilon\|^{2\sigma(1-\theta_1)}_{L^q(I_\Lambda, L^r)} \|u^\varepsilon\|^{2\sigma\theta_1}_{L^\infty(I_\Lambda, L^{2\sigma+2})}$$
$$\lesssim \varepsilon^{-2\sigma(1-\theta_1)/q} \|u^\varepsilon\|^{2\sigma\theta_1}_{L^\infty(I_\Lambda, L^{2\sigma+2})}.$$

Since $u^\varepsilon = v^\varepsilon + w^\varepsilon + r^\varepsilon$, we have,

$$\|u^\varepsilon\|^{2\sigma\theta_1}_{L^\infty(I_\Lambda, L^{2\sigma+2})} \lesssim \|w^\varepsilon\|^{2\sigma\theta_1}_{L^\infty(I_\Lambda, L^{2\sigma+2})} + \|v^\varepsilon\|^{2\sigma\theta_1}_{L^\infty(I_\Lambda, L^{2\sigma+2})} + \|r^\varepsilon\|^{2\sigma\theta_1}_{L^\infty(I_\Lambda, L^{2\sigma+2})}.$$

By the Gagliardo-Nirenberg inequality (see Lemma 2.4), we obtain as in point **1)**,

$$\|w^\varepsilon(t)\|_{L^{2\sigma+2}} \lesssim \varepsilon^{-\delta(2\sigma+2)} \|w^\varepsilon\|^{1-\delta(2\sigma+2)}_{L^2} \|\varepsilon \nabla_x w^\varepsilon\|^{\delta(2\sigma+2)}_{L^2},$$
$$\lesssim \varepsilon^{-\delta(2\sigma+2)} \|w^\varepsilon\|^{1-\delta(2\sigma+2)}_{L^2},$$

by use of the boundedness of $\|w^\varepsilon\|_{H^1_\varepsilon}$. Moreover by assumption,

$$\|r^\varepsilon\|_{L^\infty(I_\Lambda, L^{2\sigma+2})} \xrightarrow[\varepsilon \to 0]{} 0.$$

On the other hand, since $\{|t - t_j| < \Lambda\varepsilon\} \cap I_\Lambda = \emptyset$ for all $j \in \{1, ..., J\}$, we have, from Lemma 5.1,

$$\|v^\varepsilon\|^{2\sigma\theta_1}_{L^\infty(I_\Lambda, L^{2\sigma+2})} \lesssim \left(\frac{1}{\Lambda\varepsilon + \varepsilon}\right)^{n\sigma\frac{2\sigma\theta_1}{2\sigma+2}}.$$

Since

$$2\sigma(1-\theta_1)/q + n\sigma\, 2\sigma\theta_1/(2\sigma+2) = 1 - n\sigma + 2/\underline{q},$$

Eq. (5.3) reads

$$\|f^\varepsilon\|_{L^{\underline{q}'}(I_\Lambda, L^{\underline{r}'})} \lesssim \varepsilon^{-n\sigma+1+2/\underline{q}} \left(\Lambda^{-n\sigma\frac{2\sigma\theta_1}{2\sigma+2}} + o(1) + \|w^\varepsilon\|^{\gamma}_{L^\infty(I_\Lambda, L^2)}\right) \|w^\varepsilon\|_{L^q(I_\Lambda, L^2)},$$

with $\gamma = 2\sigma\theta_1 (1 - \delta(2\sigma+2)) > 0$. Hence **2)** of Lemma 5.2. $\square$

**The source term**: We prove the following lemma.

**Lemma 5.3.** *Let $T > 0$ be such that*

$$\limsup_{\varepsilon \to 0} \sup_{0 \le t \le T} \varepsilon^{n\sigma} \|r^\varepsilon(t)\|^{2\sigma+2}_{L^{2\sigma+2}} = 0.$$

*Then we have*

$$\|s^\varepsilon\|_{L^1_T(L^2)} = o(\varepsilon^{1-n\sigma}).$$

*Proof of Lemma 5.3.* Decompose the source term $s^\varepsilon$ as $s^\varepsilon = s^\varepsilon_1 + s^\varepsilon_2$ with

$$s^\varepsilon_1 := |v^\varepsilon|^{2\sigma} v^\varepsilon - \sum_{j=1}^{J} |v^\varepsilon_j|^{2\sigma} v^\varepsilon_j.$$

Let us first estimate $s^\varepsilon_2$: we have $|s^\varepsilon_2| \lesssim |v^\varepsilon|^{2\sigma}|r^\varepsilon| + |r^\varepsilon|^{2\sigma}|v^\varepsilon|$, hence by Hölder's inequality, we have

$$\|s^\varepsilon_2\|_{L^1_T(L^2)} \lesssim \|r^\varepsilon\|^{2\sigma}_{L^{2\sigma+1}_T(L^{4\sigma+2})} \|v^\varepsilon\|_{L^{2\sigma+1}_T(L^{4\sigma+2})} + \|v^\varepsilon\|^{2\sigma}_{L^{2\sigma+1}_T(L^{4\sigma+2})} \|r^\varepsilon\|_{L^{2\sigma+1}_T(L^{4\sigma+2})}.$$

By Lemma 3.3, we have

$$\|v^\varepsilon\|_{L^{2\sigma+1}_T(L^{4\sigma+2})} \lesssim \varepsilon^{-\frac{n\sigma\theta_2}{2\sigma+2}} \varepsilon^{-(1-\theta_2)/q_1}.$$

By the smallness assumption on $r^\varepsilon$, we have also

$$\|r^\varepsilon\|_{L^{2\sigma+1}_T(L^{4\sigma+2})} \le o(1) \varepsilon^{-\frac{n\sigma\theta_2}{2\sigma+2}} \varepsilon^{-(1-\theta_2)/q_1}.$$

The conclusion is as in Section 3, and we get

$$\|s^\varepsilon_2\|_{L^1_T(L^2)} \le o(\varepsilon^{1-n\sigma}).$$

Now let us consider the term $s_1^\varepsilon$. We have the pointwise estimate

$$|s_1^\varepsilon| \lesssim \sum_{j=1}^{J} |v_j^\varepsilon|^{2\sigma} \sum_{k \neq j} |v_k^\varepsilon|.$$

This is proved by induction on $J$, as soon as it is known for $J = 2$, in which case it is a consequence of the following general estimate: let $g(z) := |z|^{2\sigma}z$, then for every $(z_1, z_2) \in \mathbb{C}^2$,

$$|g(z_1 + z_2) - g(z_1) - g(z_2)| \lesssim |z_1|^{2\sigma}|z_2| + |z_2|^{2\sigma}|z_1|.$$

We first reduce the estimate to profile interactions, which is then treated by orthogonality arguments. It is proved in [Car00] that for every $\ell \in [1, J]$,

$$v_\ell^\varepsilon(t,x) = \frac{1}{\varepsilon^{n/2}} \Psi_\ell\left(\frac{t-t_\ell}{\varepsilon}, \frac{x-x_\ell}{\varepsilon}\right) + \delta_\ell^\varepsilon(t,x) =: \psi_\ell^\varepsilon(t,x) + \delta_\ell^\varepsilon(t,x),$$

where $\Psi_\ell \in C(\mathbb{R}, \Sigma) \cap L^q(\mathbb{R}, W^{1,r})$ for any admissible $(q,r)$, and

$$(5.4) \quad \forall r \in \left[2, \frac{2n}{n-2}\right[, \quad \|\delta_\ell^\varepsilon(t)\|_{L^r} \leq \frac{o(1)}{(\varepsilon + |t-t_\ell|)^{\delta(r)}}, \quad \text{uniformly in} \quad t \in \mathbb{R}.$$

It follows that

$$\|s_1^\varepsilon\|_{L_T^1(L^2)} \lesssim \sum_{j=1}^{J} \sum_{k \neq j} \left\||\psi_j^\varepsilon|^{2\sigma} \psi_k^\varepsilon\right\|_{L_T^1(L^2)} + \rho^\varepsilon,$$

where $\rho^\varepsilon$ is a linear combination of terms of the type $\left\|(|\psi_j^\varepsilon|^{2\sigma} + |\delta_\ell^\varepsilon|^{2\sigma})\delta_k^\varepsilon\right\|_{L_T^1(L^2)}$.

Let us start by showing that $\rho^\varepsilon$ is estimated by $o(\varepsilon^{1-n\sigma})$. We have indeed

$$\left\|(|\psi_j^\varepsilon|^{2\sigma} + |\delta_\ell^\varepsilon|^{2\sigma})\delta_k^\varepsilon\right\|_{L_T^1(L^2)} \leq \left(\|\psi_j^\varepsilon\|_{L_T^{2\sigma+1}(L^{4\sigma+2})}^{2\sigma} + \|\delta_\ell^\varepsilon\|_{L_T^{2\sigma+1}(L^{4\sigma+2})}^{2\sigma}\right)$$
$$\times \|\delta_k^\varepsilon\|_{L_T^{2\sigma+1}(L^{4\sigma+2})}$$

and (5.4) yields the expected result.

Finally for $j \neq k$, an obvious change of variables yields

$$\left\||\psi_j^\varepsilon|^{2\sigma} \psi_k^\varepsilon\right\|_{L_T^1(L^2)} \leq$$

$$\leq \varepsilon^{1-n\sigma} \int_{\mathbb{R}} \left(\int_{\mathbb{R}^n} \left(|\Psi_j(s,y)|^{2\sigma} \left|\Psi_k\left(s + \frac{t_j - t_k}{\varepsilon}, y + \frac{x_j - x_k}{\varepsilon}\right)\right|\right)^2 dy\right)^{1/2} ds.$$

By density one can suppose that the $\Psi_\ell$'s are compactly supported in time and space. The fact that

$$\left|\frac{t_j - t_k}{\varepsilon}\right| + \left|\frac{x_j - x_k}{\varepsilon}\right| \xrightarrow[\varepsilon \to 0]{} +\infty$$

shows that the above integral goes to zero as $\varepsilon$ goes to zero; the result follows. $\square$

*Proof of Theorem 1.4.* We now have all the tools to complete the proof of Theorem 1.4. The idea is that on intervals of the form $I_\Lambda$, we can use the linearizability condition, provided that $\Lambda$ is sufficiently large. This is quantified by the second point of Lemma 5.2. The complementary of such intervals can be split into a finite number of intervals of the form $I_\eta$ (once $\Lambda$ is fixed). Choosing first $\Lambda$ sufficiently large, then $\eta$ sufficiently small, Theorem 1.4 stems from Lemma 5.3.

Assume that the first focusing time is $t_1$. Let $0 \leq t \leq t_1 - \Lambda\varepsilon$, for $\Lambda$ to be fixed later. From (5.2) applied with $(q_1, r_1) = (q_2, r_2) = (\underline{q}, \underline{r})$, and $(q_3, r_3) = (\infty, 2)$, we have, using the second point of Lemma 5.2,

$$\|w^\varepsilon\|_{L^q_t(L^r)} \leq C_1\left(\varepsilon^{n\sigma-1-2/q}\|f^\varepsilon\|_{L^{q'}_t(L^{r'})} + \varepsilon^{n\sigma-1-1/q}\|s^\varepsilon\|_{L^1_t(L^2)}\right)$$

$$\leq C_2\left(\delta + \|w^\varepsilon\|^\gamma_{L^\infty_t(L^2)}\right)\|w^\varepsilon\|_{L^q_t(L^r)} + C_1\varepsilon^{n\sigma-1-1/q}\|s^\varepsilon\|_{L^1_t(L^2)},$$

where $\delta > 0$ is to be fixed, and the constants $C_1$, $C_2$ are universal. We first fix $\delta = 1/(4C_2)$. Lemma 5.2 yields an associate $\Lambda_0$, and we fix $\Lambda = \Lambda_0$. As in the proof of Th. 1.1, we can absorb the first term of the right hand side so long as

$$(5.5) \qquad \|w^\varepsilon\|^\gamma_{L^\infty_t(L^2)} \leq \frac{1}{4C_2}.$$

So long as (5.5) holds and $t \leq t_1 - \Lambda\varepsilon$, we deduce

$$\|w^\varepsilon\|_{L^q_t(L^r)} \leq 2C_1\varepsilon^{n\sigma-1-1/q}\|s^\varepsilon\|_{L^1_t(L^2)},$$

and like in the proof of Th. 1.1, we infer, along with Lemma 5.3, that (5.5) holds up to time $t = t_1 - \Lambda\varepsilon$, and with $I_\Lambda = [0, t_1 - \Lambda\varepsilon]$, we have

$$(5.6) \qquad \|w^\varepsilon\|_{L^\infty(I_\Lambda, L^2)} = o(1).$$

Now we have to analyze the crossing of the small time interval $[t_1 - \Lambda\varepsilon, t_1 + \Lambda\varepsilon]$. The idea is to decompose that interval into a finite number (of the order $\Lambda/\eta$) of intervals $I_\eta$ where $\eta$ is fixed so that we can repeat a similar absorption argument. Indeed, with $I_\eta = [t_1 - \Lambda\varepsilon, t_1 - \Lambda\varepsilon + \eta\varepsilon]$, we have, using Lemma 2.2

$$\begin{aligned}\|w^\varepsilon\|_{L^q(I_\eta, L^r)} &\leq C_1\varepsilon^{n\sigma-1-2/q}\|f^\varepsilon\|_{L^{q'}(I_\eta, L^{r'})} \\ &+ C_1\varepsilon^{n\sigma-1-1/q}\|s^\varepsilon\|_{L^1(I_\eta, L^2)} \\ &+ C_1\varepsilon^{-\frac{1}{q}}\|w^\varepsilon(t_1 - \Lambda\varepsilon)\|_{L^2}.\end{aligned}$$

Choosing $\delta = 1/(2C_1)$ in the first part of Lemma 5.2 fixes the value of $\eta = \eta_0$, and we have

$$\|w^\varepsilon\|_{L^q(I_\eta, L^r)} \leq 2C_1(\varepsilon^{n\sigma-1-1/q}\|s^\varepsilon\|_{L^1(I_\eta, L^2)} + \varepsilon^{-\frac{1}{q}}\|w^\varepsilon\|_{L^\infty(I_\Lambda, L^2)}).$$

We also deduce, applying Strichartz inequality once more, as in Section 3,

$$\|w^\varepsilon\|_{L^\infty(I_\eta, L^2)} \leq C(\varepsilon^{n\sigma-1}\|s^\varepsilon\|_{L^1(I_\eta, L^2)} + \|w^\varepsilon\|_{L^\infty(I_\Lambda, L^2)}).$$

which implies, using (5.6), that

$$\|w^\varepsilon\|_{L^\infty(I_\eta, L^2)} \leq C\varepsilon^{n\sigma-1}\|s^\varepsilon\|_{L^1(I_\eta, L^2)} + o(1).$$

Repeating this argument in order to cover the whole interval $[t_1 - \Lambda\varepsilon, t_1 + \Lambda\varepsilon]$ by intervals $I_\eta$, we finally come up with

$$\|w^\varepsilon\|_{L^\infty(t_1-\Lambda\varepsilon, t_1+\Lambda\varepsilon; L^2)} \leq Ce^{c\Lambda/\eta}\left(\varepsilon^{n\sigma-1}\|s^\varepsilon\|_{L^1(t_1-\Lambda\varepsilon, t_1+\Lambda\varepsilon; L^2)} + o(1)\right) = o(1),$$

from Lemma 5.3. At time $t = t_1 + \Lambda\varepsilon$, we are thus reduced to the same situation as at time $t = 0$, and we can repeat the same operations. Notice that the values of $\Lambda$ and $\eta$ are fixed once and for all at the first step, and Condition (5.5) will hold on the whole interval $[0, T]$, provided that $\varepsilon$ is sufficiently small. $\square$

*Remark.* Under the assumptions of Theorem 1.4, one cannot obtain a result in the space $L^\infty([0, T], H^1_\varepsilon)$: indeed, if we differentiate the equation on $w^\varepsilon$ and then follow the same method as above, all the terms can be estimated along the same lines except for one, which is of the type $\varepsilon^{n\sigma}|v^\varepsilon|^{2\sigma}\varepsilon\nabla r^\varepsilon$. Far from focusing times, the linearizability condition makes that term small, but the problem comes from the times $t_j$, near which $\varepsilon^{n\sigma}|v^\varepsilon|^{2\sigma}$ is not small. Under our assumptions, $\varepsilon\nabla r^\varepsilon$ has no reason to be small, in any sense. On the other hand, the asymptotics holds in $L^\infty([0, T], H^1_\varepsilon)$ if one supposes that $\varepsilon\nabla r^\varepsilon$ is linearizable (which is the case for instance if $\varepsilon\nabla r^\varepsilon_0$ goes to zero in $L^2$).

## Appendix A. Nonlinear Schrödinger equation with harmonic potential

In this appendix, we transpose the previous results to the case of nonlinear Schrödinger equations with harmonic potential. We replace (1.6) with

$$
(A.1) \quad \begin{cases} i\varepsilon \partial_t u^\varepsilon + \frac{1}{2}\varepsilon^2 \Delta u^\varepsilon = \frac{|x|^2}{2} u^\varepsilon + \varepsilon^{n\sigma} |u^\varepsilon|^{2\sigma} u^\varepsilon, & (t,x) \in \mathbb{R}_+ \times \mathbb{R}^n, \\ u^\varepsilon_{|t=0} = u^\varepsilon_0. \end{cases}
$$

We still suppose that $\sigma > 2/n$, and $\sigma < 2/(n-2)$ if $n \geq 3$. The motivation in this study relies in the fact that such equations are currently used to model Bose-Einstein condensation (see e.g. [CT99]). In [KNSQ00], it is proposed that the above equation models Bose-Einstein condensation in space dimensions one ($n = 1$, with $\sigma = 2$) and two ($n = 2$, with $\sigma = 1$), if $\varepsilon = \hbar$, the Planck constant. In space dimension three, the nonlinearity $\hbar^2 |u^\hbar|^2 u^\hbar$ is usually considered, and does not fit our scaling. Notice that (A.1) meets the usual model when $n \leq 2$, and $\sigma = 2/n$, which is precisely the borderline case for which we do not know whether all our results still hold. Nevertheless, some results remain, and we believe that the information provided by the case $\sigma > 2/n$ is interesting.

The corresponding linear solution satisfies

$$
(A.2) \quad \begin{cases} i\varepsilon \partial_t v^\varepsilon + \frac{1}{2}\varepsilon^2 \Delta v^\varepsilon = \frac{|x|^2}{2} v^\varepsilon, \\ v^\varepsilon_{|t=0} = u^\varepsilon_0. \end{cases}
$$

It is natural to assume not only that $u^\varepsilon_0 \in H^1_\varepsilon(\mathbb{R}^n)$, but also that $|x| u^\varepsilon_0 \in L^2(\mathbb{R}^n)$, so that $u^\varepsilon_0$ belongs to the domain of $\sqrt{-\varepsilon^2 \Delta + |x|^2}$. We denote

$$
\Sigma_\varepsilon := \left\{ f^\varepsilon \in H^1(\mathbb{R}^n) \; ; \; \sup_{0 < \varepsilon \leq 1} \left( \|f^\varepsilon\|_{L^2} + \|\varepsilon \nabla_x f^\varepsilon\|_{L^2} + \|x f^\varepsilon\|_{L^2} \right) < \infty \right\}.
$$

As before, we also suppose that there is no focusing at time 0,

$$
\limsup_{\varepsilon \to 0} \varepsilon^{n\sigma} \|u^\varepsilon_0\|^{2\sigma+2}_{L^{2\sigma+2}} = 0.
$$

The expression of $v^\varepsilon$ is given by Mehler's formula (see e.g. [FH65]),

$$
v^\varepsilon(t,x) = \frac{1}{(2i\pi\varepsilon \sin t)^{n/2}} \int_{\mathbb{R}^n} e^{\frac{i}{\varepsilon \sin t}\left(\frac{|x|^2+|y|^2}{2}\cos t - x \cdot y\right)} u^\varepsilon_0(y) dy =: U^\varepsilon(t) u^\varepsilon_0(x).
$$

In particular, local dispersion estimates hold for the group $U^\varepsilon$ and the same Strichartz estimates as in the case with no potential hold, with a constant depending on the size of the time interval considered (see for instance [Caz93], Remark 3.4.4). Following [Car03], introduce

$$
A^\varepsilon(t) = x \sin t - i\varepsilon \cos t \nabla_x, \quad B^\varepsilon(t) = x \cos t + i\varepsilon \sin t \nabla_x.
$$

These two operators, usual in a linear context (they are the quantization of impulse and momentum), commute with the operator

$$
i\varepsilon \partial_t + \frac{1}{2}\varepsilon^2 \Delta - \frac{|x|^2}{2},
$$

act on the nonlinearity $|u^\varepsilon|^{2\sigma} u^\varepsilon$ like derivatives, and satisfy the pointwise identity

$$
(A.3) \quad |A^\varepsilon(t)\varphi(x)|^2 + |B^\varepsilon(t)\varphi(x)|^2 = |x\varphi(x)|^2 + |\varepsilon \nabla_x \varphi(x)|^2.
$$

Using the above remarks, we have the following preliminary results.

**Lemma A.1.** *Assume that $u_0^\varepsilon \in \Sigma_\varepsilon$. Let $I$ be a bounded time interval. The functions $u^\varepsilon$ and $v^\varepsilon$ are the unique solutions of (A.1) and (A.2) in $C(\mathbb{R}; \Sigma_\varepsilon)$ respectively, and satisfy the following properties.*

(1) $u^\varepsilon, v^\varepsilon \in C(\mathbb{R}; \Sigma_\varepsilon)$, *with the following conservations.*
  - *Mass:* $\|u^\varepsilon(t)\|_{L^2} = \|v^\varepsilon(t)\|_{L^2} = \|u_0^\varepsilon\|_{L^2}$.
  - *Linear energy:*
$$E_0^\varepsilon(t) := \frac{1}{2}\|\varepsilon\nabla_x v^\varepsilon(t)\|_{L^2}^2 + \frac{1}{2}\|xv^\varepsilon(t)\|_{L^2}^2 = E_0^\varepsilon(0).$$
  - *Nonlinear energy:*
$$E^\varepsilon(t) := \frac{1}{2}\|\varepsilon\nabla_x u^\varepsilon(t)\|_{L^2}^2 + \frac{1}{2}\|xu^\varepsilon(t)\|_{L^2}^2 + \frac{\varepsilon^{n\sigma}}{\sigma+1}\|u^\varepsilon(t)\|_{L^{2\sigma+2}}^{2\sigma+2} = E^\varepsilon(0).$$

(2) *For any admissible pair $(q,r)$, there exists $C_r(I)$ independent of $\varepsilon$ such that*

(A.4) $\quad \|v^\varepsilon\|_{L^q(I;L^r)} + \|A^\varepsilon v^\varepsilon\|_{L^q(I;L^r)} + \|B^\varepsilon v^\varepsilon\|_{L^q(I;L^r)} \leq C_r(I)\varepsilon^{-1/q}.$

(3) *For any admissible pair $(q_1, r_1)$ and $(q_2, r_2)$, there exists $C_{r_1,r_2}(I)$ independent of $\varepsilon$ such that*

$$\varepsilon^{\frac{1}{q_1}+\frac{1}{q_2}}\left\|\int_{I\cap\{s\leq t\}} U_0^\varepsilon(t-s)F(s)ds\right\|_{L^{q_1}(I;L^{r_1})} \leq C_{r_1,r_2}(I)\|F\|_{L^{q_2'}(I;L^{r_2'})}.$$

We did not state Strichartz estimates for $u^\varepsilon$ in (A.4). First, we noticed at the end of Sect. 3 that they are not really needed, for knowing such estimates for $v^\varepsilon$ is enough. On the other hand we could get these estimates for $u^\varepsilon$ thanks to the following change of unknown (see [Car02]). Define $\widetilde{u}^\varepsilon$ by

(A.5) $\quad \widetilde{u}^\varepsilon(t,x) = \frac{1}{(1+t^2)^{n/4}} e^{i\frac{t}{2\varepsilon(1+t^2)}|x|^2} u^\varepsilon\left(\arctan t, \frac{x}{\sqrt{1+t^2}}\right).$

Then $\widetilde{u}^\varepsilon$ solves

$$i\varepsilon\partial_t\widetilde{u}^\varepsilon + \frac{1}{2}\varepsilon^2\Delta\widetilde{u}^\varepsilon = \varepsilon^{n\sigma}(1+t^2)^{\frac{n\sigma}{2}-1}|\widetilde{u}^\varepsilon|^{2\sigma}\widetilde{u}^\varepsilon.$$

Since we have to consider bounded values for $t$, Strichartz estimates for $u^\varepsilon$ are a consequence of Lemma 2.2, along with the remark that

$$A^\varepsilon\left(t+\frac{\pi}{2}\right) = B^\varepsilon(t).$$

Using either one of the above two arguments, we have the first result,

**Theorem A.2.** *Let $T > 0$. The following properties are equivalent,*
(1) *The function $v^\varepsilon$ is an approximation of $u^\varepsilon$ on the time interval $[0,T]$,*

$$\sup_{0\leq t\leq T}\left(\|u^\varepsilon(t) - v^\varepsilon(t)\|_{L^2} + \|\varepsilon\nabla_x u^\varepsilon(t) - \varepsilon\nabla_x v^\varepsilon(t)\|_{L^2} + \|xu^\varepsilon(t) - xv^\varepsilon(t)\|_{L^2}\right) \xrightarrow[\varepsilon\to 0]{} 0.$$

(2) *The function $v^\varepsilon$ satisfies*

(A.6) $\quad \limsup_{\varepsilon\to 0}\sup_{0\leq t\leq T}\varepsilon^{n\sigma}\|v^\varepsilon(t)\|_{L^{2\sigma+2}}^{2\sigma+2} = 0.$

This result is proved the same way as Theorem 1.1. Notice that from (A.3), (1) is equivalent to

$$\sup_{0\leq t\leq T}\left(\|u^\varepsilon(t) - v^\varepsilon(t)\|_{L^2} + \|A^\varepsilon(t)(u^\varepsilon - v^\varepsilon)\|_{L^2} + \|B^\varepsilon(t)(u^\varepsilon - v^\varepsilon)\|_{L^2}\right) \xrightarrow[\varepsilon\to 0]{} 0.$$

Since $A^\varepsilon$ and $B^\varepsilon$ commute with the linear equation and act on the nonlinearity like derivatives, each of them plays the same role as $\varepsilon\nabla_x$ in Theorem 1.1.

*Remark.* Notice that as in Theorem 1.1, the proof that (1) ⇒ (2) does not rely on the assumption $\sigma > 2/n$, since it is a consequence of the conservation laws. In particular, this result could be used to study Bose-Einstein condensation in space dimensions one and two.

To state the analog of Theorem 1.2, notice that (A.5) turns (A.2) into (1.7).

**Theorem A.3.** *Let $T > 0$ and assume that (1.9) is not satisfied. Then:*
*(1) There exists $0 < T_0 \leq \pi$ such that*
$$\limsup_{\varepsilon \to 0} \sup_{0 \leq t \leq T_0} \varepsilon^{n\sigma} \|v^\varepsilon(t)\|_{L^{2\sigma+2}}^{2\sigma+2} > 0.$$

*(2) Up to the extraction of a subsequence, there exist an orthogonal family $(t_j^\varepsilon, x_j^\varepsilon)_{j \in \mathbb{N}}$ in $]0, \pi[ \times \mathbb{R}^n$, a family $(\Psi_\ell^\varepsilon)_{\ell \in \mathbb{N}}$, bounded in $H_\varepsilon^1(\mathbb{R}^n)$, and a (nonempty) family $(\varphi_j)_{j \in \mathbb{N}}$ bounded in $\mathcal{F}(H^1)(\mathbb{R}^n)$, such that for all $\ell \in \mathbb{N}$ and all $x \in \mathbb{R}^n$,*
$$u_0^\varepsilon(x) = \Psi_\ell^\varepsilon(x) + r_\ell^\varepsilon(x), \quad \text{with } \limsup_{\varepsilon \to 0} \|U^\varepsilon(t) r_\ell^\varepsilon\|_{L^\infty_{\text{loc}}(\mathbb{R}_+, L^{2\sigma+2}_\varepsilon)} \xrightarrow[\ell \to \infty]{} 0,$$

*and for every $\ell \in \mathbb{N}$, the following asymptotics holds in $L^2(\mathbb{R}^n)$, as $\varepsilon \to 0$,*
$$\Psi_\ell^\varepsilon(x) = \sum_{j=0}^\ell \frac{1}{(\sin t_j^\varepsilon)^{\frac{n}{2}}} \varphi_j\left(\frac{x - x_j^\varepsilon \cos t_j^\varepsilon}{\sin t_j^\varepsilon}\right) e^{-i\frac{|x|^2 + |x_j^\varepsilon|^2}{2\varepsilon} \cot t_j^\varepsilon + i \frac{x \cdot x_j^\varepsilon}{2\varepsilon \sin t_j^\varepsilon}} + o(1).$$

*Moreover, we have $\limsup_{\varepsilon \to 0} \frac{t_j^\varepsilon}{\varepsilon} = +\infty$ for all $j \in \mathbb{N}$, and there is an integer $j \in \mathbb{N}$ such that $t_j^\varepsilon \in [0, T_0]$.*

*Remark.* In particular, if the nonlinear term remains negligible up to time $T = \pi$, then it is *always* negligible. Notice that the first point of the theorem is actually a consequence of the profile decomposition stated in the second point. This phenomenon can be compared with a result stated in [Car02], where the special case $\sigma = 2/n$ is considered. Let $u$ be the solution of the initial value problem (notice that the nonlinearity is attractive, and not repulsive as in the present paper),
$$\begin{cases} i\partial_t u + \frac{1}{2}\Delta u = \frac{|x|^2}{2} u - |u|^{4/n} u, & (t, x) \in \mathbb{R}_+ \times \mathbb{R}^n, \\ u_{|t=0} = u_0. \end{cases}$$

If the mass of $u_0$ is critical (that is, if its $L^2$-norm is that of the solitary wave associated to (1.2)), then blow up in finite time may occur. But if the solution has not collapsed up to time $t = \pi$, then it never blows up.

*Remark.* Such initial data still make the nonlinear term relevant at time $t = t_j^\varepsilon$ (concentration at the point $x_j^\varepsilon$) in the case $\sigma = 2/n$, in particular for Bose-Einstein condensation. The open question is: are they the only ones?

**Corollary A.4.** *Assume that for every Wigner measure $\mu_0$ associated with the data $u_0^\varepsilon$ and for every $y \in \mathbb{R}^n$,*
$$\mu_0 \perp \delta(x - y) \otimes d\xi.$$

*Assume moreover that there exists $0 < T \leq \pi$ such for every $y \in \mathbb{R}^n$ and for every $t \in ]0, T[$,*
$$\mu_0 \perp dx \otimes \delta(\xi - x \cot t - y).$$

*Then*
$$\limsup_{\varepsilon \to 0} \sup_{0 \leq t \leq T} \varepsilon^{n\sigma} \|v^\varepsilon(t)\|_{L^{2\sigma+2}}^{2\sigma+2} = 0.$$

*In other words, $u^\varepsilon$ is linearizable on $[0, T]$. Moreover, if one can take $T = \pi$, then $u^\varepsilon$ is linearizable on $\mathbb{R}_+$.*

We finally analyze the case where the initial data $u_0^\varepsilon$ has the form displayed in Th. A.3,

$$\text{(A.7)} \quad u_0^\varepsilon(x) = \sum_{j=1}^J \frac{1}{(\sin t_j)^{\frac{n}{2}}} f_j\left(\frac{x - x_j \cos t_j}{\sin t_j}\right) e^{-i\frac{|x|^2 + |x_j|^2}{2\varepsilon} \cot t_j + i\frac{x \cdot x_j}{2\varepsilon \sin t_j}} + r_0^\varepsilon(x),$$

where $f_j \in \Sigma$, $x_j \in \mathbb{R}^n$, $t_j \in ]0, \pi[$, with $(t_j, x_j) \neq (t_k, x_k)$ if $j \neq k$. We also assume that $r_0^\varepsilon$ is bounded in $\Sigma_\varepsilon$ and that its free evolution $r^\varepsilon$, defined by

$$\text{(A.8)} \quad \begin{cases} i\varepsilon \partial_t r^\varepsilon + \frac{1}{2}\varepsilon^2 \Delta r^\varepsilon = \frac{|x|^2}{2} r^\varepsilon, \\ r^\varepsilon_{|t=0} = r_0^\varepsilon, \end{cases}$$

satisfies (A.6) for some $0 < T \leq \pi$. For $1 \leq j \leq J$, we define $v_j^\varepsilon$ as the solution of the initial value problem

$$\text{(A.9)} \quad \begin{cases} i\varepsilon \partial_t v_j^\varepsilon + \frac{1}{2}\varepsilon^2 \Delta v_j^\varepsilon = \frac{|x|^2}{2} v_j^\varepsilon + \varepsilon^{n\sigma} |v_j^\varepsilon|^{2\sigma} v_j^\varepsilon, \quad (t,x) \in \mathbb{R}_+ \times \mathbb{R}^n, \\ v_{j|t=0}^\varepsilon = \frac{1}{(\sin t_j)^{\frac{n}{2}}} f_j\left(\frac{x - x_j \cos t_j}{\sin t_j}\right) e^{-i\frac{|x|^2 + |x_j|^2}{2\varepsilon} \cot t_j + i\frac{x \cdot x_j}{2\varepsilon \sin t_j}}. \end{cases}$$

The asymptotics for $v_j^\varepsilon$ is described in [Car03] for $t_j = \pi/2$. For $t_j \neq \pi/2$, it is similar and can be deduced by a time translation.

**Theorem A.5.** *Assume that $u_0^\varepsilon$ is given by (A.7). Then for any $T > 0$ such that*

$$\limsup_{\varepsilon \to 0} \sup_{0 \leq t \leq T} \varepsilon^{n\sigma} \|r^\varepsilon(t)\|_{L^{2\sigma+2}}^{2\sigma+2} = 0,$$

*the following asymptotics holds in $L^\infty(0, T; L^2)$ as $\varepsilon$ goes to zero,*

$$u^\varepsilon = \sum_{j=1}^J v_j^\varepsilon + r^\varepsilon + o(1).$$

*Remark.* All the above results could be generalized to the case of an anisotropic potential,

$$V(x) = \frac{1}{2}\left(\omega_1^2 x_1^2 + \omega_2^2 x_2^2 + \ldots + \omega_n^2 x_n^2\right),$$

with $\omega_j \geq 0$, not necessarily all equal. Mehler's formula in that case is different but analogous, hence Strichartz estimates are still available. The technical point then consists in replacing the operators $A^\varepsilon$ and $B^\varepsilon$ by their vectorial counterparts defined by

$$A_j^\varepsilon(t) = \omega_j x_j \sin(\omega_j t) - i\varepsilon \cos(\omega_j t) \partial_j, \quad B_j^\varepsilon(t) = \omega_j x_j \cos(\omega_j t) + i\varepsilon \sin(\omega_j t) \partial_j.$$

The asymptotics of the corresponding $v_j^\varepsilon$'s are discussed in some particular cases in [Car03].

MAB et UMR CNRS 5466, Université Bordeaux 1, 351 cours de la Libération, 33 405 Talence cedex, France
   *E-mail address*: carles@math.u-bordeaux.fr

Université de Cergy-Pontoise, Mathématiques, 2 avenue Adolphe Chauvin, BP 222, Pontoise, 95 302 Cergy-Pontoise cedex, France
   *E-mail address*: Clotilde.Fermanian@math.u-cergy.fr

Centre de Mathématiques, UMR CNRS 7640, École Polytechnique, 91 128 Palaiseau cedex, France
   *E-mail address*: Isabelle.Gallagher@math.polytechnique.fr